\documentclass[11pt]{amsart}
   \newlength{\extramargin}
\usepackage{amsfonts, amsmath,amscd, amssymb, latexsym, mathrsfs, stmaryrd, verbatim, wasysym }
\usepackage[all, 2cell, ps]{xy}
\usepackage{graphicx}
\DeclareGraphicsExtensions{.ps}
\newtheorem{theorem}{Theorem}
\newtheorem{corollary}[theorem]{Corollary}

\newtheorem{lemma}[theorem]{Lemma}
\newtheorem{proposition}[theorem]{Proposition}

\numberwithin{equation}{section}
\numberwithin{theorem}{section}

 \newcommand{\Fun}{
 \makebox{
 $\text{\scriptsize{I}} \hspace{-.3 mm} \text{\scriptsize{I}}$
 \hspace{-.2 in} \raisebox{.05 in}{$\smile$}
 \hspace{-.23 in} \raisebox{-.045 in}{$\frown$}}}

 \newcommand{\ITFun}{
 \makebox{
 $\text{\scriptsize{I}} \hspace{-.35 mm} \text{\scriptsize{I}}$
 \hspace{-.19 in} \raisebox{.05 in}{$\smile$}
 \hspace{-.25 in} \raisebox{-.045 in}{$\frown$}}}

 \newcommand{\SupFun}{
 \makebox{
 $\text{\scriptsize{I}} \hspace{-.3 mm} \text{\scriptsize{I}}$
 \hspace{-.214 in} \raisebox{.047 in}{$\smile$}
 \hspace{-.249 in} \raisebox{-.055 in}{$\frown$}}}

 \newcommand{\SupITFun}{
 \makebox{
 $\text{\scriptsize{I}} \hspace{-.35 mm} \text{\scriptsize{I}}$
 \hspace{-.205 in} \raisebox{.05 in}{$\smile$}
 \hspace{-.28 in} \raisebox{-.05 in}{$\frown$}}}

 \newcommand{\FunForm}{
 \makebox{$\mathcal{I} \hspace{-1.25 mm} \mathcal{I}$}}
 \newcommand{\SFun}{
 \makebox{$\mathrm{I}  \hspace{-.16 mm} \mathrm{I}$}}
 \newcommand{\vSFun}{
 \makebox{$\mathbb{I}  \hspace{-.16 mm} \mathbb{I}$}}

\DeclareMathAlphabet{\mathpzc}{OT1}{pzc}{m}{it}
\DeclareMathOperator{\cl}{\mathpzc{cl}}
\DeclareMathOperator{\hook}{\lrcorner}

\DeclareMathOperator*{\FP}{\text{FP}}

\DeclareMathOperator*{\btimes}{\times}
\DeclareMathOperator{\Vol}{\text{Volume}}
\DeclareMathOperator{\dvol}{\text{dvol}}

\DeclareMathOperator{\Ind}{\text{Ind}}

\renewcommand{\mod}{\text{ }\mathrm{mod}\text{ }}
\renewcommand{\bar}[1]{\overline{#1}}

\newcommand{\Rp}{\mathbb{R}^+}
\newcommand{\Hom}{\text{Hom}}

\newcommand{\script}[1]{\textsc{#1}}
\newcommand{\df}[1]{\mathfrak{#1}}
\newcommand{\curly}[1]{\mathcal{#1}}
\newcommand{\bhs}[1]{\mathfrak{B_{#1}}}

\newcommand{\lrpar}[1]{\left( #1 \right)}
\newcommand{\lrspar}[1]{\left[ #1 \right]}

\title[Gauss-Bonnet Theorem]{A renormalized index theorem for some complete asymptotically regular metrics: the Gauss-Bonnet theorem}
\author{Pierre Albin}
\address{Department of Mathematics, Massachusetts Institute of Technology}
\email{pierre@math.mit.edu} 

\newcommand\datver[1]{\def\datverp%
 {\par\boxed{\boxed{\text{Version: #1; Run: \today}}}}}
\datver{1.0; Revised: \today}

\begin{document}
%%%%%%%%%%%%%%%%%%%%%%%%%%%%%%%%%%%%%%%%%%%%%%%%%%%%%%%%%%%%%%%%%%%%%%%%%%%%%%%%%%%%%%%%%%%%%%%%%%%%%%%%%%%%%%%%%%%%%%%%%%%%

\date\datverp

\begin{abstract}
The Gauss-Bonnet Theorem is studied for edge metrics as a renormalized index theorem. These metrics include the Poincar\'e-Einstein metrics of the AdS/CFT correspondence. Renormalization is used to make sense of the curvature integral and the dimensions of the $L^2$-cohomology spaces as well as to carry out the heat equation proof of the index theorem.
For conformally compact metrics even mod $x^m$, the finite time supertrace of the heat kernel on conformally compact manifolds is shown to renormalize independently of the choice of special boundary defining function.
\end{abstract}

\maketitle

%%%%%%%%%%%%%%%%%%%%%%%%%%%%%%%%%%%%%%%%%%%%%%%%%%%%%%%%%%%%%%%%%%%%%%%%%%%%%%%%%%%%%%%%%%%%%%%%%%%%%%%%%%%%%%%%%%%%%%%%%%%%
\section{Introduction}
%%%%%%%%%%%%%%%%%%%%%%%%%%%%%%%%%%%%%%%%%%%%%%%%%%%%%%%%%%%%%%%%%%%%%%%%%%%%%%%%%%%%%%%%%%%%%%%%%%%%%%%%%%%%%%%%%%%%%%%%%%%%

In \cite{APS}, Atiyah, Patodi, and Singer extended the Atiyah Singer index 
theorem from closed manifolds to manifolds with boundary. If $\eth$ is 
a Dirac-type operator, then
\[ \int AS - \frac{1}{2}\eta(\partial M) = h + \Ind(\eth). \]
Whereas on a closed manifold the index is computed solely as the integral of the Atiyah-Singer ``curvature integrand'', $AS$, in the presence of a boundary 
one must include ``extended'' 
solutions, denoted by $h$, and 
the eta invariant of the boundary metric, $\eta(\partial M)$. A key observation 
in \cite{APS} shows the equivalence of the APS \textit{spectral} boundary 
condition and the $L^2$ condition for ``asymptotically cylindrical'' non-compact 
manifolds.

Not all of the index theorems on non-compact spaces involve Fredholm operators having an {\em a priori} well-defined index, and one makes use of an averaging procedure (such as the group action in Atiyah's $\Gamma$-index theorem) or a regular exhausion (in \cite{Roe:Ind}) to interpret the dimension of the infinite dimensional kernel and cokernel.
%These theorems are referred to as type II in \cite{Roe:Ind}. 
The more traditional index theorems %are known as type I and
 include, along with the Atiyah-Patodi-Singer theorem above, its reworking by Melrose 
in \cite{APS Book}, extension to manifolds with corners in 
\cite{Corners} and \cite{Salomonsen}, extension to locally symmetric spaces 
in \cite{Stern}, and Carron's theory of Dirac operators which are 
``non-parabolic at infinity'' \cite{Carron:Ind}, among many others.

%In particular, Dirac-type operators on hyperbolic space (or, more generally, conformally compact manifolds, as defined below) do not fit well into either of these categories.

An interesting case not fitting well into either of these categories is that of Dirac-type operators on
 conformally compact manifolds, such as hyperbolic space and Poincar\'e-Einstein spaces.
The problem is that the these operators (e.g., the Gauss-Bonnet and signature operators) have 
infinite dimensional kernel and/or cokernel. 
In the present work, we work on more general ``edge" type metrics, and interpret the index by a renormalization discussed below.

The relevance of Poincar\'e-Einstein, or PE, manifolds in conformal geometry and the AdS/CFT correspondence in physics has fueled a lot of recent investigations. 
Among these, the study of the Gauss-Bonnet theorem was carried out in four dimensions in \cite{Anderson}, for constant sectional curvature PE manifolds in \cite{Epstein} and in general even dimensions in \cite{Chang-Qing-Yang} and \cite{My Preprint}. In contrast to previous treatments of this theme, here we wish to consider it as a (renormalized) index theorem. The complications arise because the $L^2$-cohomology is infinite dimensional (in degree $\frac{1}{2} \dim M$).

%%%%%%%%%%%%%%%%%%%%%%%%%%%%%%%%%%%%%%%%%%%%%%%%%%%%%%%%%%%%%%%%%%%%%%%%%%%%%%%%%%%%%%%%%%%%%%%%%%%%%%%%%%%%%%%%%%%%%%%%%%%%
\subsection{Summary and statement of results} $ $\newline
%%%%%%%%%%%%%%%%%%%%%%%%%%%%%%%%%%%%%%%%%%%%%%%%%%%%%%%%%%%%%%%%%%%%%%%%%%%%%%%%%%%%%%%%%%%%%%%%%%%%%%%%%%%%%%%%%%%%%%%%%%%%

A convenient setting is afforded by complete metrics on the interior, $M$, of a compact manifold with boundary, $\bar{M}$. This facilitates a discussion of asymptotic regularity and is general enough to allow the study of, for instance, asymptotically cylindrical, asymptotically hyperbolic, and asymptotically locally Euclidean as well as edge metrics.

One can make use of a boundary defining function (or bdf: a function $x \in C^\infty\lrpar{\bar{M}}$ with a simple zero at $\partial M$ and positive otherwise) to identify a neighborhood of $\partial M$ with a product
\begin{equation} \label{ProdNhd}
	[0,1)_x \times \partial M.
\end{equation}
In such a neighborhood, asymptotically hyperbolic metrics can be written
\begin{equation*}
	\frac{dx^2}{x^2} + h_x.
\end{equation*}
It is often useful to study the product of an asymptotically hyperbolic metric and a closed manifold. Edge metrics are asymptotically modeled by twisted versions of these products.

For these metrics, the boundary is the total space of a fibration,
\begin{equation*}
	F - \partial M \xrightarrow{\phi} B,
\end{equation*}
and the vector fields of bounded point-wise length are precisely those tangent to the fibers of the fibration at the boundary. 
A complete edge metric (a product edge metric in the nomenclature of \cite{Vaillant}) is one that in a neighborhood like \eqref{ProdNhd} takes the form
\begin{equation} \label{ProdEdge}
	\frac{dx^2}{x^2} + \frac{\phi^*\lrpar{g_B}}{x^2} + g_F,
\end{equation}
where $g_B$ is a metric on the base and $g_F$ is a symmetric two-tensor that restricts to a metric on each fiber.
We will deal exclusively with complete edge metrics and henceforth refer to them simply as edge metrics.
When the dimension of the base, $b$, is zero, edge metrics are asymptotically cylindrical and when the dimension of the fiber, $f$, is zero, they are asymptotically hyperbolic.

The elliptic theory of edge metrics is described in \cite{Mazzeo:Edge}, while \cite{APS Book} contains both the elliptic theory and the index theorem for asymptotically cylindrical metrics. Our proof of the index theorem for edge metrics is based on this proof of the Atiyah-Patodi-Singer Index Theorem.
We begin by reviewing the elliptic edge calculus in section $\mathcal{x}$\ref{sec:elliptic}. In sections $\mathcal{x}$\ref{sec:HeatSpace} and $\mathcal{x}$\ref{sec:Comp} we set up the ``edge heat calculus" with bundle coefficients, $E$, and determine its composition properties.

\begin{proposition}\label{Intro1}
There is a bi-filtered algebra of integral operators, $\Psi^{k,\ell}_{e, Heat}(M,E)$ with composition rule
\begin{equation*}
	\Psi^{k_1,\ell_1}_{e, Heat}(M,E) \circ \Psi^{k_2,\ell_2}_{e, Heat}(M,E)
	\subset \Psi^{k_1+k_2, \ell_1 + \ell_2}_{e, Heat}(M,E)
\end{equation*}
and ``symbol maps" consistent with the filtrations, known as normal operators, such that
\begin{equation}\begin{split}
      &0 \to \Psi^{k-1, \ell}_{e, Heat} \to \Psi^{k, \ell}_{e, Heat} \xrightarrow{N_{00,2}} 
           \rho^{\ell}_{11,0}  \dot{\curly{C}}^{\infty}_{11,0}
            \lrpar{\bhs{00,2}, E\downharpoonright_{\bhs{00,2}}} \to 0 \\
      &0 \to \Psi^{k, \ell-1}_{e, Heat} \to \Psi^{k, \ell}_{e,Heat} \xrightarrow{N_{11,0}} 
            \rho^{k}_{00,2}  \dot{\curly{C}}^{\infty}_{00,2}
            \lrpar{\bhs{11,0}, E\downharpoonright_{\bhs{11,0}}} \to 0
\end{split} \end{equation}
are short exact sequences.
\end{proposition}
Here $\bhs{11,0}$ and $\bhs{00,2}$ are two boundary hypersurfaces of the ``heat space" of  $\mathcal{x}$\ref{sec:HeatSpace}, and $\dot{\curly{C}}^{\infty}_{X}$ denotes smooth functions vanishing to infinite order at all boundary hypersurfaces other than $X$.

Given these normal operators, we construct a solution to the heat equation in much the same way as one uses the symbol map of pseudo-differential operators to construct parametrices of elliptic operators on closed manifolds.

\begin{theorem}\label{Intro2}
The heat kernel associated to the Laplacian of an edge metric is an element of $\curly{H} \in \Psi^{2,0}_{e,Heat}(M,E)$.
\end{theorem}

In section $\mathcal{x}$\ref{sec:IndexTheorem}, we present a proof of the index theorem.
The heat equation proof of the index theorem on a closed manifold boils down to showing that the supertrace of the heat kernel is independent of time and comparing its values as time approaches zero and infinity. The heat kernel in Theorem \ref{Intro2} is not trace-class nor is the Laplacian Fredholm, so in $\mathcal{x}$\ref{sec:RInt} we describe a scheme for renormalizing integrals, traces and dimensions.
The Laplacian of an edge metric (with $b>0$) has a spectral gap at zero, so as time approaches infinity the renormalized supertrace converges to the renormalized index of the de Rham operator. However a generalized Dirac operator might not have closed image, e.g., this is the case for the Dirac operator on Hyperbolic space \cite{Bunke}.
As time approaches zero, the local index theorem asserts the convergence of the pointwise supertrace of the heat kernel to the integrand of the Atiyah-Singer index theorem, $AS$.
Unlike for closed manifolds, the renormalized supertrace of the heat kernel does depend on time, so the index theorem contains an extra term, ${}^R\eta$, which in the asymptotically cylindrical case is the eta invariant of \cite{APS}.
We show that ${}^R\eta$ localizes to a neighborhood of the boundary, and vanishes for the de Rham operator. The full analysis of this invariant will be carried out in \cite{Me-Rafe} through a Getzler rescaling of the bundle at the boundary, along the lines of \cite[Chapter 8]{APS Book}.

\begin{theorem} \label{Intro3}
If $M$ is an edge manifold (with $\dim B>0$) and $\eth^E$ is a generalized Dirac operator,
\begin{equation*}
	{}^R\int_M AS + {}^R\eta = \lim_{t\to \infty} {}^R\mathrm{Str}\lrpar{e^{-\lrpar{\eth^E}^2}}.
\end{equation*}
In the special case of the Gauss-Bonnet operator this limit exists and equals the renormalized index, ${}^R\mathrm{Ind}\lrpar{\eth_{GB}}$. The renormalized eta invariant vanishes and hence
\begin{equation*}
	{}^R\int \mathrm{Pff} = {}^R\mathrm{Ind}\lrpar{\eth_{GB}}.
\end{equation*}
\end{theorem}

Notice that this theorem applies in particular to the conformally compact metrics.
In the asymptotically cylindrical case ($b=0$), a full analysis can be found in \cite{APS Book}. The index of $\eth^E$ is always finite though the image is closed only when the induced Dirac operator over the boundary has no null space. The corresponding index theorem is that of \cite{APS}.

We focus on the Gauss-Bonnet theorem in section $\mathcal{x}$\ref{sec:GBThm}. 
The main task is to specifically interpret the two sides of this theorem. For the ${}^R\int \mathrm{Pff}$, we prove a `soft' index theorem. 
By considering the Chern-Gauss-Bonnet formula on the truncated manifold with boundary $\{x\geq \varepsilon\}$:
\begin{equation*}
	\int_{x\geq\varepsilon} \mathrm{Pff} + \int_{x = \varepsilon} \SupFun 
	= \chi \lrpar{ \{ x \geq \varepsilon \} },
\end{equation*}
and letting $\varepsilon \to 0$, we obtain the following result.

\begin{theorem} \label{Intro4}
For an edge metric,
\begin{equation*}
	{}^R\int_{M} \mathrm{Pff} + \FP_{\varepsilon=0}\int_{x = \varepsilon} \SupITFun 
	= \chi \lrpar{ M }.
\end{equation*}
Thus the renormalized index differs from the Euler characteristic by a boundary integral,
\begin{equation*}
	{}^R\mathrm{Ind}\lrpar{\eth_{GB}} = \chi(M) - \FP_{\varepsilon=0}\int_{x = \varepsilon} \SupITFun .
\end{equation*}
\end{theorem}

Einstein edge metrics with trivial fibre ($f=0$), or Poincar\'e-Einstein metrics, are of interest to the physics community because of their role in the AdS/CFT correspondence. We showed in \cite{My Preprint} that the boundary integral in Theorem \ref{Intro4} vanishes in this context.
One interesting aspect of renormalization for these metrics is the existence of a distinguished subset of bdfs. Many invariants have the same renormalization independently of the choice the bdf from among these ``special" bdfs. In \cite{My Preprint}, this was seen to be the case for all scalar Riemannian invariants. We were also able to show that $\FP_{\varepsilon=0}\int_{x = \varepsilon} \SupFun$ so that,
in conjunction with the fact that on any conformally compact manifold, \cite{Mazzeo:Hodge}
\begin{equation*}\begin{split}
	{}^R\mathrm{Ind}\lrpar{\eth_{GB}} &=
         \sum_{k<\frac{m}{2}} (-1)^k \dim H^k(M,\partial M)
         +(-1)^{m/2}  \lrpar{ {}^R\dim \mathcal{H}_{L^2}^{m/2}}
         +\sum_{k>\frac{m}{2}} (-1)^k \dim H^k(M)\\
         &= 2\sum_{k<\frac{m}{2}} (-1)^k \dim H^k(M,\partial M)
	+ (-1)^{m/2}  \lrpar{ {}^R\dim \mathcal{H}_{L^2}^{m/2}},
\end{split} \end{equation*}
the following is true.

\begin{corollary} 
For Poincar\'e-Einstein metrics (or any conformally compact metric even below $x^m$ cf. section $\mathcal{x}$\ref{sec:ConfComp}),
\begin{equation*}
	{}^R\int \mathrm{Pff} = \chi(M).
\end{equation*}
In particular, this implies that
\begin{equation*}
	(-1)^{m/2}  \lrpar{ {}^R\dim \mathcal{H}_{L^2}^{m/2}}
	= \sum_{k \geq \frac{m}{2}} (-1)^k \dim H^k(M,\partial M)
	 - \sum_{k<\frac{m}{2}} (-1)^k \dim H^k(M,\partial M).
\end{equation*}
\end{corollary}

In section $\mathcal{x}$\ref{sec:ConfComp}, we show that the same conclusion holds for the trace of the heat kernel at any fixed finite time.

\begin{theorem} \label{Intro5}
If $g$ is a conformally compact metric that is even mod $x^m$, the renormalized trace of the heat kernel of the Hodge Laplcian acting on forms, obtained after choosing a special bdf is independent of the choice of special bdf.
\end{theorem}

The proof proceeds by constructing an even subcalculus and showing that the Laplacian and then the resolvent are elements of this subcalculus. The result for the heat kernel follows from its functional calculus expression in terms of the resolvent. The same construction could be used for other functions of the Laplacian. 
Colin Guillarmou has independently singled out the same class of even operators in his forthcoming study of a generalized Krein spectral shift function.

%%%%%%%%%%%%%%%%%%%%%%%%%%%%%%%%%%%%%%%%%%%%%%%%%%%%%%%%%%%%%%%%%%%%%%%%%%%%%%%%%%%%%%%%%%%%%%%%%%%%%%%%%%%%%%%%%%%%%%%%%%%%
\subsection{Extension to other metrics} $ $\newline
%%%%%%%%%%%%%%%%%%%%%%%%%%%%%%%%%%%%%%%%%%%%%%%%%%%%%%%%%%%%%%%%%%%%%%%%%%%%%%%%%%%%%%%%%%%%%%%%%%%%%%%%%%%%%%%%%%%%%%%%%%%%

Edge metrics are a particular example of ``boundary fibration structures" \cite{Kyoto}. Other examples include metrics that 
on a neighborhood like \eqref{ProdNhd} have the form
\begin{equation*}
	i) \frac{dx^2}{x^2} + h_x, \phantom{xxx}
	ii) \frac{dx^2}{x^4} + h_x, \phantom{xxx}
	iii) \frac{dx^2}{x^4} + \frac{h_x}{x^2}, \phantom{xxx}
	iv) \frac{dx^2}{x^4} + \frac{\phi^*g_B}{x^2} + g_F, \phantom{xxx}
	v) \frac{dx^2}{x^2} + \phi^*g_B + x^2g_F.
\end{equation*}
These model a wide variety of asymptotic geometries.
Metrics $(i)$ and $(ii)$ correspond to different compactifications of cylindrical ends and are studied in the $b$-calculus and the $cusp$-calculus, respectively. The radial compactification of $\mathbb{R}^n$ yields a metric of type $(iii)$ as do other asymptotically conical manifolds. These are the subject of the scattering-calculus. Finally, the last two metrics correspond to fibered boundaries and fibered cusps respectively.

One can proceed further and consider iterated fibrations at the boundary, with metrics of the form
\begin{equation*}
	\frac{dx^2}{x^{a_0}} + \frac{g_1}{x^{a_1}} + \cdots + \frac{g_k}{x^{a_k}}.
\end{equation*}
Those with $a_0 = \max \{a_i\}$, and all $a_i \geq 0$, can be considered as iterated edge structures, and will be refered to as MICE (Metrics with Iterated Complete Edge structures).
The construction of the heat calculus and the proof of the composition rules (Proposition \ref{Intro1}) extend to cover MICE, as explained at the end of section \ref{sec:Comp}. 
Assuming that the model problems can be solved, one should be able then use the heat kernel to prove a renormalized index theorem as before:
\begin{equation*}
	{}^R\int AS - \frac{1}{2}{}^R\eta = \lim_{t \to \infty} {}^R\mathrm{Str}\lrpar{e^{-t\eth^2}}.
\end{equation*}
For the Gauss-Bonnet operator, the proofs of the vanishing of the renormalized eta invariant and the limiting or ``soft" index formula (Theorem \ref{Intro4}) hold for MICE.
In section \ref{sec:TrivFib}, we describe the boundary integral when each of the fibrations have trivial base or fiber. In particular, we have the following result.
\begin{theorem} \label{Intro6}
For asymptotically cylindrical metrics, or more generally those with $a_i=0$ for $i>0$,
\begin{equation*}
	{}^R\int \mathrm{Pff} = \chi(M).
\end{equation*}
For scattering metrics (type $(iii)$ above),
\begin{equation*}
	\int \mathrm{Pff} + P\lrpar{\partial M, h_0} = \chi(M),
\end{equation*}
where $P\lrpar{\partial M, h_0}$ is a linear combination of the Weyl volume of tubes invariants of the boundary with the metric $h_0$.
\end{theorem}

There is an $L^2$-index theorem for asymptotically locally Euclidean metrics \cite{CarronL2},
\begin{equation*}
	\int \mathrm{Pff} = \chi_{L^2}(M).
\end{equation*}
The $L^2$-Euler characteristic is the alternating sum of the (finite) dimensions of the spaces of harmonic forms. It is known that this is a topological invariant \cite{HHM}.

\begin{corollary}
For scattering metrics,
the difference between the $L^2$-Euler characteristic and the topological Euler characteristic is given by a linear combination of the Weyl volume of tubes invariants,
\begin{equation*}
	\chi(M) - \chi_{L^2}(M) = P\lrpar{\partial M, h_0}.
\end{equation*}
\end{corollary}

\subsection{Acknowledgements} 
This work forms part of my thesis. I am very grateful to my advisor, Rafe Mazzeo, for
 his guidance and inspiration.
Throughout this work, I received support from his NSF grant DMS-0204730.
I am also grateful to Richard Melrose for fruitful conversations on this project.
I would like to thank Colin Guillarmou for uncovering an error in an earlier draft of the final section, and for useful conversations about renormalization.

%%%%%%%%%%%%%%%%%%%%%%%%%%%%%%%%%%%%%%%%%%%%%%%%%%%%%%%%%%%%%%%%%%%%%%%%%%%%%%%%%%%%%%%%%%%%%%%%%%%%%%%%%%%%%%%%%%%%%%%%%%%%
\section{The Edge Calculus} \label{sec:elliptic}
%%%%%%%%%%%%%%%%%%%%%%%%%%%%%%%%%%%%%%%%%%%%%%%%%%%%%%%%%%%%%%%%%%%%%%%%%%%%%%%%%%%%%%%%%%%%%%%%%%%%%%%%%%%%%%%%%%%%%%%%%%%%

We start by recalling the edge calculus from \cite{Mazzeo:Edge}.
Let $M$ be the interior of an $m$ dimensional manifold with boundary $\overline{M}$. Assume that the boundary $\partial M$ is the total space of a fibration
\[ F \to \partial M \to B ,\]
and denote by $\curly{V}_e \subset TM$ the vector bundles that are tangent to the fibers of the fibration.
Denote the dimensions of $M$, $F$, and $B$ by $m$, $f$, and $b$ respectively.

It is instructive to use local coordinates to describe $\curly{V}_e$. 
Let $\{ \partial_{y_1}, \ldots, \partial_{y_b} \}$ be local coordinate vector fields for $B$ and $\{ \partial_{z_1}, \ldots \partial_{z_f} \}$ local coordinate vector fields for $F$. Assume that $x$ is a boundary defining function, or bdf. That is, $x$ is a smooth positive function on $\overline{M}$ with a simple zero at $\partial M$. A full set of coordinates for $TM$ near the boundary is given by $\{ \partial_x, \partial_{y_1}, \ldots, \partial_{y_b}, \partial_{z_1}, \ldots, \partial_{z_f} \} =: \{X_i\}$, thus $\curly{V}_e$ has as local spanning set
\begin{equation}\label{edgeVF}
	\{ x\partial_x, x\partial_{y_1}, \ldots, x\partial_{y_b}, \partial_{z_1}, \ldots, \partial_{z_f} \}.
\end{equation}
As discussed in \cite{Mazzeo:Edge}, these vector fields are a basis of a bundle, known as the edge-tangent bundle, ${}^eTM$. It is isomorphic to the tangent bundle, though not canonically so. The natural inclusion
\begin{equation} \label{EdgeInc}
	{}^eTM \hookrightarrow TM
\end{equation}
is an isomorphism over the interior of $M$, but degenerates at the boundary.

An edge metric is a symmetric two tensor in the dual bundle to ${}^eTM$, ${}^eT^*M$.
In a product neighborhood like \eqref{ProdNhd}, an ``exact" edge metric (cf. \cite{APS Book}, \cite {Vaillant}) has the form
\begin{equation*}
	\frac{dx^2}{x^2} + \frac{h}{x^2} + g_F
\end{equation*}
where $h \downharpoonright_{\partial M} = \phi^*\lrpar{g_B}$ for some metric on the base $g_B$ and 
$g_F$ is a symmetric two tensor in the usual cotangent bundle $T^*M$. Thus, 
 these are asymptotically like the product edge metrics, \eqref{ProdEdge}, and at the boundary ${}^eTM$ decomposes orthogonally into 
 \begin{equation} \label{BdyDecomp}
	{}^eTM\downharpoonright_{\partial M} =
	 \left< x\partial_x, x\partial_{y_1}, \ldots, x\partial_{y_b} \right>
	 \oplus \left< \partial_{z_1}, \ldots, \partial_{z_f} \right>
\end{equation}
where the bundles on the right are the kernel and the image respectively of the map \eqref{EdgeInc} at the boundary.

The enveloping algebra of ${}^eTM$ is the space is edge differential operators. Thus, in local coordinates, an edge differential operator is a polynomial in the vector fields \eqref{edgeVF},
\begin{equation}\label{EdgeDiff}
	L=\sum_{i+|J|+|K| \leq \ell} 
	a_{i,J,K}\lrpar{x,y,z} \lrpar{x\partial_x}^i \lrpar{x\partial_y}^J \lrpar{\partial_z}^K.
\end{equation}
Analogously to the closed case, the highest order terms transform like elements in ${}^eT^*M$. As a homogeneous polynomial on ${}^eT^*M$, the symbol is invariantly defined with a local expression corresponding to \eqref{EdgeDiff} given by
\begin{equation*}
	{}^e\sigma\lrpar{L} \lrpar{\xi,\eta,\zeta}= \sum_{i+|J|+|K| = \ell} 
	a_{i,J,K}\lrpar{x,y,z} \lrpar{\xi}^i \lrpar{\eta}^J \lrpar{\zeta}^K.
\end{equation*}
An edge differential operator is elliptic if this symbol is invertible away from the zero section. The Laplacian of an edge metric is elliptic as an edge differential operator.

As usual, the search for a parametrix for elliptic differential operators leads to considering more general, pseudo-differential, operators. These operate via integration against an integral kernel with a conormal singularity along the diagonal, as for closed manifolds. However, these parametrices will generally have additional singularities where the diagonal hits the boundary. The point of view adopted in \cite{Mazzeo:Edge} is that the integral kernel should be considered a push-forward of a simpler distribution on a more complicated space which covers $M^2$ and is known as the edge-stretch product or the edge double-space.

To construct the edge double-space, we observe that the additional singularities of parametrices of edge elliptic operators occur along approaches to the submanifold of the corner,
\begin{equation*}
	\partial M \btimes_B \partial M = \{ (\zeta, \zeta') \in \partial M^2 : \phi(\zeta)=\phi(\zeta') \},
\end{equation*}
known as the fibered product of $\partial M$ with itself over $B$. To resolve these singularities we introduce polar coordinates around this submanifold. In the language of \cite{APS Book}, this corresponds to blowing-up the submanifold, and can be viewed as replacing $\partial M \btimes_B \partial M$ with its inward-pointing spherical normal bundle to obtain the edge double-space
\begin{equation*}
	M^2_e = \lrspar{ M^2, \partial M \btimes_B \partial M }.
\end{equation*}
The new boundary face produced by this blow-up is known as the (edge) front-face and is denoted $\bhs{11}$. The other boundary faces, $\bhs{10}$ and $\bhs{01}$, are (induced by) $\partial M \times M$ and $M \times \partial M$ respectively. 

The blow-up construction furnishes us with a blow-down map
\begin{equation*}
	M^2_e \xrightarrow{\beta_e} M^2,
\end{equation*}
which collapses $\bhs{11}$. This is an example of a ``$b$-map", since the pull-back of any bdf is a product of bdfs and positive functions.

The kernels of edge pseudo-differential operators are polyhomogeneous conormal distributions on $M^2_e$, as we now describe. Denote bdfs for all boundary hypersurfaces of a space $X$ by $x_i$, then using multi-index notation, a distribution $u$ is polyhomogeneous conormal if
\begin{equation}\label{phg}
	u \in \curly{A}^{*}_{phg}\lrpar{X} \implies
	u \sim \sum_{\Re s_j \to \infty} \sum_{p=0}^{p_j} 
	a_{j,p}(x,y) x^{s_j} \lrpar{\log x}^p, a_{j,p} \in \curly{C}^{\infty}.
\end{equation}
A finer space can be described by fixing the set of exponents $\{s_j, p \}$ that are allowed to occur. We require of such a discrete set, $E \subset \mathbb{C} \times \mathbb{N}_0$, that
\begin{equation*}\begin{split}
	i)& \lrpar{ s_j, p_j } \in E, |\lrpar{ s_j, p_j }| \to \infty \implies \Re\lrpar{s_j} \to \infty \\
	ii)& \lrpar{s,p} \in E \implies \lrpar{s + k, p - \ell} \in E, 
	\text{ for any } k,\ell \in \mathbb{N}, \ell \leq p,
\end{split}\end{equation*}
in which case we refer to $E$ as an ``index set" or a ``smooth index set". If $\curly{E}$ is a collection of index sets, one per boundary hypersurface, the space $\curly{A}^{\curly{E}}_{phg}$ is defined by restricting the expansions in \eqref{phg} to those having exponents from $\curly{E}$.
A distribution has a conormal singularity at an interior ``p-submanifold" (see \cite[App. A]{Mazzeo:Edge}), $Y$, if its transverse Fourier transform is a symbol. The order of the singularity is determined by the order of the symbol and the space of such functions of order $m$ is denoted by $I^m\lrpar{X,Y}$. Finally, given $\curly{E}$ and $m$, we can define the space of edge pseudo-differential operators $\Psi_e^{\curly{E}, m}(M;\Omega^{1/2})$ as those having integral kernel in $\curly{A}^{\curly{E}}_{phg}I^m\lrpar{M^2_e,\mathrm{diag}_e;\rho_{11}^{\frac{b+1}{2}}\Omega^{1/2}}$, where
\begin{equation*}
	\mathrm{diag}_e := \bar{\beta_e^{-1}\lrpar{\mathrm{diag} %
		\setminus \partial M \btimes_B \partial M}},
\end{equation*}
and $\rho_{11}$ is a bdf for $\bhs{11}$.

The reason for using half-densities is so that it makes sense to compose operators. The composition can be analyzed geometrically. Starting with the formula for the kernel of $A \circ B$,
\begin{equation*}
      \curly{K}_{A \circ B} \lrpar{\zeta, \zeta''} 
         = \int_{\zeta'} \curly{K}_A \lrpar{ \zeta, \zeta'} \curly{K}_B \lrpar{ \zeta', \zeta'' },
\end{equation*}
note that in terms of the maps
\begin{equation*}
      \xymatrix{ &  \lrpar{\zeta, \zeta', \zeta''} %
       \ar@{|->}[dl]^{\beta_{LM}} \ar@{|->}[d]^{\beta_{LR}} \ar@{|->}[dr]^{\beta_{MR}}  &  \\
         \lrpar{\zeta, \zeta'} &  \lrpar{\zeta, \zeta''} &  \lrpar{\zeta', \zeta''} } 
\end{equation*}
it can be written
\begin{equation}\label{comp}
      \curly{K}_{A \circ B} = \lrpar{\beta_{LR}}_*
             \lrpar{ \beta_{LM}^*\curly{K}_A \cdot \beta_{MR}^* \curly{K}_B } .
\end{equation}
This analysis is carried out in \cite{Mazzeo:Edge} by constructing a ``triple space", 
$M^3_e$, 
with nice maps ($b$-fibrations) down to the ``left", ``right" and ``center" double spaces,
\begin{equation}\label{Data}
      \xymatrix{ &  M^3_e \ar[dl]^{\beta_{LM}} \ar[d]^{\beta_{LR}} \ar[dr]^{\beta_{MR}} &  \\
         M^2_e &  M^2_e &  M^2_e } .
\end{equation}
$b$-fibrations are nice in this context because they preserve the space of polyhomogeneous conormal distributions (for a discussion of these, see e.g. \cite[App. A]{Mazzeo:Edge}, \cite{Grieser}, or \cite{Corners} ). The triple space is constructed by blowing up the three copies of $\partial M \btimes_B \partial M$ to be found in $M^3$, which we denote $S_{LM}$, $S_{LR}$, and $S_{MR}$, after having blown-up their intersection, $S_{LMR}$. Thus,
\begin{equation*}
	M^3_e = \lrspar{ M^3; S_{LMR}; S_{LM} \cup S_{LR} \cup S_{MR} }.
\end{equation*}

\begin{figure}[htpb]
      \centering
      \includegraphics[bb=0in 2.45in 8.6in 8.7in, keepaspectratio, height=2in,clip]{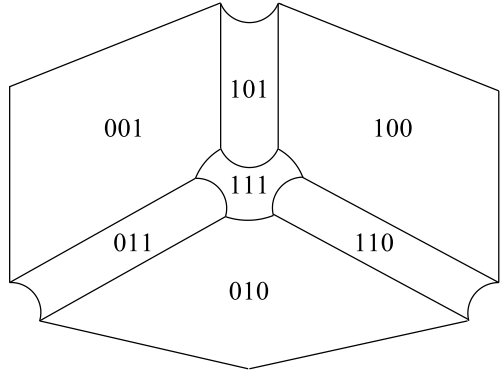}
      \caption{A diagram of the triple space, $M^3_e$.}
      \label{Triple space}
\end{figure}

The composition result shown in \cite{Mazzeo:Edge} is as follows. If $A \in \Psi^{m,\curly{E}}_e$, $B \in \Psi^{m', \curly{F}}_e$ and 
\begin{equation*}
	\Re\lrpar{E_{01}} + \Re\lrpar{F_{10}} > -1,
\end{equation*}
then $A \circ B$ is defined and is an element of $\Psi^{m+m',\curly{G}}_e$ with
\begin{equation*}\begin{split}
	G_{10} &= \lrpar{E_{11} + F_{10}} \bar{\cup} E_{10} \\
	G_{01} &= \lrpar{E_{01} + F_{11}} \bar{\cup} F_{01} \\
	G_{11} &= \lrpar{E_{11} + F_{11}} \bar{\cup} \lrpar{E_{10} + F_{01} + k +1}.
\end{split}\end{equation*}
Here the operation $\bar{\cup}$ or extended union of the index sets is defined by
\begin{equation*}
	E \bar{\cup} F = E \cup F \cup
		\{ \lrpar{z,p+p'+1} : \lrpar{z,p} \in E, \lrpar{z,p'}\in F \}.
\end{equation*}

In \cite{Mazzeo:Edge} it is shown that the calculus above is large enough to contain parametrices of elliptic edge operators, when these exist.
These are constructed with the aid of a second ``symbol map", the normal operator, $N(A)$, obtained by restricting the kernel of $A$ to the edge front-face, $\bhs{11}$. An operator is Fredholm if both its classical symbol and its normal operator are invertible. Since the Laplacian is a natural operator, its normal operator is in turn the Laplacian of the model metric. Hence, for an exact edge metric, by \eqref{BdyDecomp}
\begin{equation}\label{NLap}
	N_e^2(\Delta) = \Delta_{\mathbb{H}^{b}}\Delta_{F}.
\end{equation}
This has implications for the essential spectrum of the Laplacian. It is easy to see that this only depends on the behavior near the boundary. Indeed, multiplication by a function $\chi$ supported in the interior of the manifold is a compact operator so
\begin{equation*}
	\sigma_{ess}\lrpar{\Delta} = \sigma_{ess}\lrpar{\Delta \circ \curly{M}_{1-\chi}}.
\end{equation*}
More to the point, since $\Delta -\lambda$ is elliptic, it will be Fredholm precisely when $N_e(\Delta- \lambda)$ is invertible.

In particular, since the Laplacian of a conformally flat metric has a spectral gap at zero \cite{Mazzeo:Hodge}, so will the Laplacian of an exact edge metric (with $\dim B >0$). On the other hand, since the Dirac operator of hyperbolic space does not have a spectral gap at zero \cite{Bunke}, neither will that of a general edge metric.

%%%%%%%%%%%%%%%%%%%%%%%%%%%%%%%%%%%%%%%%%%%%%%%%%%%%%%%%%%%%%%%%%%%%%%%%%%%%%%%%%%%%%%%%%%%%%%%%%%%%%%%%%%%%%%%%%%%%%%%%%%%%
\section{The Heat Space} \label{sec:HeatSpace}
%%%%%%%%%%%%%%%%%%%%%%%%%%%%%%%%%%%%%%%%%%%%%%%%%%%%%%%%%%%%%%%%%%%%%%%%%%%%%%%%%%%%%%%%%%%%%%%%%%%%%%%%%%%%%%%%%%%%%%%%%%%%

It is convenient to keep in mind the heat kernel of Euclidean space,
\begin{equation*}
      h_0(\zeta, \zeta', t) = \frac{1}{(2\pi t)^{\frac{n}{2}}} \exp\lrpar{ -\frac{| \zeta - \zeta'|^2}{4t} } .
\end{equation*}
Away from $\{ \zeta = \zeta' \}$, the heat kernel vanishes to infinite order as $t \to 0$, but along the diagonal at $\{t=0\}$ it blows up. To understand the behavior along different directions of approach to this submanifold $\{ \zeta = \zeta', t=0\}$, the invariance of $t^{-1}| \zeta - \zeta'|^2$ with respect to $(\zeta, \zeta', t) \mapsto (\lambda \zeta, \lambda \zeta', \lambda^2 t)$ suggests performing a parabolic blow-up.

This heuristic is true also for the heat kernel of an edge metric. 
We define the heat space corresponding to $M^2_e$ by performing a parabolic blow-up, in the direction of $dt$, of the blown-up diagonal at time zero. Symbolically,
\begin{equation*}
     HM_e = \lrspar{ M^2_e \times \mathbb{R}^+ ; \mathrm{diag}_e \times \{0\}, \left< dt \right> } .
\end{equation*}

\begin{figure}[htpb]
      \centering
      \includegraphics[bb=0in 2.75in 8.5in 8.5in, keepaspectratio, height=2in,clip]{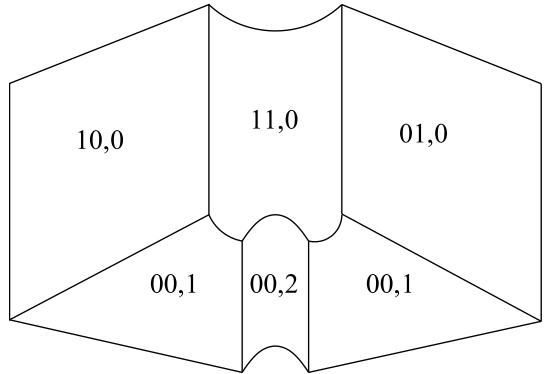}
      \caption{A diagram of heat space, $HM_e$.}
      \label{HeatSpace}
\end{figure}

The construction of $HM_e$ furnishes us with a natural blow-down map,
\begin{equation*}
      HM_e \xrightarrow{\beta_H} M^2_e \times \mathbb{R}^+ ,
\end{equation*}
which collapses the boundary face created by the blow-up, 
$\bhs{00,2}$,
 the ``temporal front face". 
The boundary faces of $HM_e$ are easily described. The faces $\bhs{10}$, $\bhs{01}$, $\bhs{11}$ from $M^2_e$ give rise to boundary hypersurfaces which we denote by $\bhs{10,0}$, $\bhs{01,0}$, and $\bhs{11,0}$ respectively. 
The ``temporal boundary face", $\bhs{00,1}$, consists of that part of $\{ t= 0\}$ that was not blown-up, formally
\begin{equation*}
      \bhs{00,1} = \overline{ \beta_H^{-1}\lrpar{ \lrpar{ M^2 \setminus \mathrm{diag}_e} \times \{0\} } }.
\end{equation*}
We shall use $\rho_{*}$ to denote boundary defining functions of bdfs. Thus, for instance, $\rho_{00,2}$ is a positive function vanishing simply at $\bhs{00,2}$.

The effect of the blow-up is most clearly seen in terms of local coordinates.
We describe projective coordinates, valid away from $\{ x'=0 \}$.
In the interior, we can use the usual coordinates
\[ \lrpar{ \lrpar{ x,y,z }, \lrpar{ x',y',z'}, t} = \lrpar{ \zeta, \zeta', t } .\]
These same coordinates, with $t$ replaced by $t^{1/2}$ work away from $\bhs{11,0}$.
Near $\bhs{11,0}$, but away from $\bhs{00,1}$, we use
\begin{equation}\label{efCoord}
       \lrpar{ \lrpar{ \frac{x}{x'}, \frac{y-y'}{x'}, z}, \lrpar{x', y', z'}, t^{1/2} }
       =: \lrpar{ \lrpar{s,u,z}, \zeta', t^{1/2} } .
\end{equation}
Finally, near both $\bhs{11,0}$ and $\bhs{00,2}$,
\[ \lrpar{ \lrpar{ \frac{x - x'}{x't^{1/2}}, \frac{y-y'}{x't^{1/2}}, \frac{z-z'}{t^{1/2}} }, \lrpar{x', y', z'}, t^{1/2} }
       =: \lrpar{ \lrpar{\curly{S},\curly{U},\curly{Z}}, \zeta',t^{1/2} } \]
are smooth coordinates.

We describe the (edge) heat calculus by specifying the regularity of its integral kernels. Define the density bundle, $KD$, by
\begin{equation}\label{KD}
      KD := \rho_{00,2}^{-\lrpar{\frac{m+3}{2}}} \rho_{11,0}^{-\lrpar{\frac{b+1}{2}}}
                 \Omega^{1/2}\lrpar{HM_e}.
\end{equation}
The kernels of the operators in the heat calculus will be elements of
\begin{equation}\label{DefK}
     \mathscr{K}^{k,\ell}(M,KD) := \rho^{k}_{00,2}\rho^{\ell}_{11,0} 
                                      \dot{\curly{C}}^{\infty}_{00,2; 11,0}\lrpar{HM_e, KD}, k>0, \ell \geq 0
\end{equation}
where $\dot{\curly{C}}^{\infty}_{00,2; 11,0}$ refers to functions vanishing to infinite order at all boundary hypersurfaces except $\bhs{00,2}$ and $\bhs{11,0}$. For the moment we will assume that  $k>0$, we deal with $k=0$ by introducing additional ``mean value" assumptions, as in \cite{DM}, at the end of this section.

These kernels define operators in a couple of different ways. Recall that the solution to the heat equation in Euclidean space is given by
\begin{equation*}
      \begin{cases}
          \lrpar{\partial_t + \Delta}f(t,x) = 0\\
          \displaystyle \lim_{t \to 0}f(x,t) = f(x)
      \end{cases}
      \implies
           f(x,t) = \int_{\mathbb{R}^n} h_0(x, x',t) f(x') dx' .
\end{equation*}
Analogously, given an element $\curly{K}_A \in \mathscr{K}^{k,\ell}$, we define an operator
\begin{equation}\label{1stOp}
      \curly{C}^{\infty}\lrpar{M;\Omega^{1/2} } \xrightarrow{A}
           \curly{C}^{-\infty}\lrpar{M \times \mathbb{R}^+;\Omega^{1/2}}
      \text{ by }  A\lrpar{f}( \zeta,t) = 
            \int_M \lrpar{\gamma_* \curly{K}_A}(\zeta, \zeta',t)f(\zeta'),
\end{equation}
where $\gamma$ denotes the full blow-down map, $HM_e \to M^2_e \times \Rp \to M^2 \times \Rp$.
We denote the space of such operators with integral kernels in $\mathscr{K}^{k,\ell}(M,KD)$ by $\Psi^{k,\ell}_{e, Heat}(M,KD)$.

In \eqref{1stOp}, $\curly{C}^{-\infty}$, the dual space of $\curly{C}^{\infty}$, is a much larger target than we need. To understand the actual regularity of this mapping, we take a closer look at the push-forward of the kernel.
First note that a partition of unity allows us to assume that $f$ has small support, and that for $f$ supported away from $\bhs{00,2}$ and $\bhs{11,0}$, $A(f)$ is clearly in $\dot{\curly{C}}^{\infty}$. So we may assume that $f$ has small support near $\bhs{00,2}$.
It will be convenient to use the coordinates
\begin{equation}\label{ConvCoord}
      \lrpar{ \zeta, \lrpar{ \frac{x-x'}{xt^{1/2}}, \frac{y-y'}{xt^{1/2}}, \frac{z-z'}{t^{1/2}} }, t^{1/2} }=:
          \lrpar{ \zeta, \lrpar{ \curly{S}', \curly{U}', \curly{Z}'}, t^{1/2}} ,
\end{equation}
for which we may use $x$ for $\rho_{11,0}$ and $t^{1/2}$ for $\rho_{00,2}$.

Recall \cite{Grieser} that push-forward of a density is the dual of pull-back, that is,
\begin{equation*}
      \int \phi \beta_*(\mu) = \int \beta^*(\phi) \mu .
\end{equation*}
Thus, if we write 
\[ \curly{K}_A = 
        \kappa_A  \frac{|d\zeta d\curly{S}' d\curly{U}' d\curly{Z}' dt^{1/2}|^{1/2}}{t^{%
        \frac{m+3}{4}-\frac{k}{2}} x^{\frac{b+1}{2}-\ell} } =
        \kappa_A  \frac{|d\zeta d\curly{S}' d\curly{U}' d\curly{Z}' dt|^{1/2}}{t^{%
        \frac{m}{4}-\frac{k}{2}+1} x^{\frac{b+1}{2}-\ell} }        ,\]
we can use that the Jacobian of 
\begin{equation*}
      \zeta'  \to \lrpar{\curly{S}', \curly{U}', \curly{Z}'}  \text{ is } x^{b+1}t^{m/2},
\end{equation*}
to write the formula in \eqref{1stOp} as
\begin{equation} \label{EffectiveAct} \begin{split}
      &\lrpar{Af}(\zeta, t) = \int \lrpar{\gamma}_*\curly{K}_A(\zeta, \zeta', t)f(\zeta') %
           = \int \kappa_A 
           \frac{|d\zeta d\curly{S}' d\curly{U}' d\curly{Z}' dt|^{1/2}}{%
                    t^{\frac{m+4}{4}-\frac{k}{2}}x^{\frac{b+1}{2}-\ell} }
           \lrpar{\gamma \circ \pi_R}^*f           \\
      &= \lrspar{ \int%
              \kappa_A\lrpar{\zeta, \curly{S}',\curly{U}',\curly{Z}', t^{1/2}}%
              f\lrpar{x-xt^{1/2}\curly{S}',y-xt^{1/2}\curly{U}',z-t^{1/2}\curly{Z}'}%
              d\curly{S}' d\curly{U}' d\curly{Z}' }
          \frac{ |d\zeta dt|^{1/2} }{ t^{1-\frac{k}{2}}x^{-\ell} },
\end{split}\end{equation}
where $\pi_R:M^2 \to M$ is the projection onto the right factor.
We end up with a much more satisfying version of \eqref{1stOp},
\begin{equation}\label{1stOpDeux}
      A \in \Psi^{k, \ell}_{e, Heat}(M,KD) \implies
      \curly{C}^{\infty}\lrpar{M;\Omega^{1/2}_M } \xrightarrow{A}
          t^{\frac{k}{2}-1}x^{\ell} 
          \curly{C}^{\infty}\lrpar{M \times \mathbb{R}^+_{1/2};\Omega^{1/2}_M |dt|^{1/2}},
\end{equation}
where the $1/2$ in $\mathbb{R}^+_{1/2}$ indicates that the functions are smooth in $t^{1/2}$ instead of $t$. 

There is another way in which these operators act. Consider the bilinear map \cite[(B.16)]{DM}, \cite[(3.22)]{DF}
\begin{equation*}
     \dot{\curly{C}}^{\infty}_0\lrpar{M \times \mathbb{R}^+; \Omega^{1/2} } \times
     \dot{\curly{C}}^{\infty}\lrpar{M \times \mathbb{R}^+;\Omega^{1/2} } \xrightarrow{\hat{*}_t}
     \dot{\curly{C}}^{\infty}\lrpar{M^2 \times \mathbb{R}^+;\Omega^{1/2} }
\end{equation*}
given by
\begin{equation*}
      \phi \hat{*}_t \psi = \int_0^{\infty} \phi\lrpar{\zeta, t+t'} \psi\lrpar{\zeta',t'} dt' .
\end{equation*}
An operator $A \in \Psi^{k,\ell}_{e, Heat}\lrpar{M,KD}$ defines a continuous linear map
\begin{equation}\label{2ndOp}
      {\curly{C}}^{-\infty}\lrpar{M \times \mathbb{R}^+; \Omega^{1/2} } \xrightarrow{\widetilde{A}}
      {\curly{C}}^{-\infty}\lrpar{M \times \mathbb{R}^+; \Omega^{1/2} }  
\end{equation}
by demanding that
\begin{equation*}
      \left< \widetilde{A}\psi , \phi \right> = \int_{HM_e} \lrpar{\curly{K}_A} \beta_H^*\lrpar{\phi \hat{*}_t \psi} .
\end{equation*}
Equivalently, 
\begin{equation}\label{2ndOpForm}
      \widetilde{A}\psi\lrpar{\zeta,t} =
          \int_{M}\int_{0}^t \lrpar{\gamma_*\curly{K}_A} \lrpar{\zeta, \zeta' ,s} 
             \psi \lrpar{\zeta', t-s} ds d\zeta' .
\end{equation}
Note that if $\psi \in \dot{C}^{\infty}\lrpar{M;\Omega^{1/2}}$, then
\begin{equation*}
      \widetilde{\psi}\lrpar{\zeta, t} := \psi(\zeta) \delta(t) \implies
      \widetilde{A}\widetilde{\psi} = A\psi ,
\end{equation*}
so solving the heat equation means finding an operator, $A$, such that
\begin{equation*}
      \begin{cases}
          \lrpar{\partial_t + \Delta}A = 0\\
          \displaystyle \lim_{t \to 0}A = \text{Id}
      \end{cases}
      \iff  \lrpar{\partial_t + \Delta } \widetilde{A} = \delta(t) \otimes \text{Id} .
\end{equation*}

Notice that the right hand side of the last equation is precisely the kernel of the identity as a convolution operator \eqref{2ndOpForm}. The kernel density bundle \eqref{KD} is defined precisely so that the identity is (so far formally) an operator of order zero.
Define the space of kernels 
\begin{equation*}
           \mathscr{K}^{0,\ell}(M,KD) := \rho^{\ell}_{11,0} 
                                      \dot{\curly{C}}^{\infty}_{\overline{00,2}; 11,0}\lrpar{HM_e, KD}, \ell \geq 0,
\end{equation*}
where the bar over $\overline{00,2}$ indicates the following mean value condition. It makes sense to restrict any kernel $\curly{K} \in \mathscr{K}^{0,\ell}(M,KD)$ to $\bhs{00,2}$ which fibers over $M$ (with fiber a compactified edge tangent bundle). We require that the integral of $\curly{K}$ along each fiber vanishes.

The point is that for a kernel in $\curly{K} \in \mathscr{K}^{0,\ell}(M,KD)$,
the integral in \eqref{2ndOp} does not necessarily converge near $\bhs{00,2}$. On the other hand, if we restrict to a neighborhood of $\bhs{00,2}$, say $\curly{U}:=[0,\varepsilon) \times {}^eTM$, then
\begin{equation*}
      \int_{\curly{U}} \curly{K} \cdot \beta_H^*\lrpar{\phi \hat{*}_t \psi}
      = \lim_{\delta \to 0} \int_{\delta}^{\varepsilon} \int_M \int_{\text{fibre}}
         F_{\curly{K},\phi, \psi} \dvol \frac{dt^{1/2}}{t^{1/2}},
\end{equation*}
exists because, by the mean value condition, the innermost integral goes to zero with $t$ (see \cite[(72)]{Vaillant}).

%We shall prove \eqref{2ndOp} later, for the moment 
An advantage of \eqref{2ndOp} over \eqref{1stOpDeux} is that these operators can be composed. We shall analyze the composition in the next section.

%%%%%%%%%%%%%%%%%%%%%%%%%%%%%%%%%%%%%%%%%%%%%%%%%%%%%%%%%%%%%%%%%%%%%%%%%%%%%%%%%%%%%%%%%%%%%%%%%%%%%%%%%%%%%%%%%%%%%%%%%%%%
\section{Composition of Operators}\label{sec:Comp}
%%%%%%%%%%%%%%%%%%%%%%%%%%%%%%%%%%%%%%%%%%%%%%%%%%%%%%%%%%%%%%%%%%%%%%%%%%%%%%%%%%%%%%%%%%%%%%%%%%%%%%%%%%%%%%%%%%%%%%%%%%%%

The composition formula for heat operators acting by convolution, as in \eqref{2ndOp}, is
\begin{equation*}
      \curly{K}_{A \circ B} \lrpar{\zeta, \zeta'', t} 
         = \int_{\zeta',t'} \curly{K}_A \lrpar{ \zeta, \zeta', t-t' } \curly{K}_B \lrpar{ \zeta', \zeta'', t' } .
\end{equation*}
To analyze this formula geometrically, note that in terms of the maps
\begin{equation*}
      \xymatrix{ &  \lrpar{\zeta, \zeta', \zeta'', t-t', t'} %
       \ar@{|->}[dl]^{\beta_{LM,t-t'}} \ar@{|->}[d]^{\beta_{LR,t}} \ar@{|->}[dr]^{\beta_{MR,t'}}  &  \\
         \lrpar{\zeta, \zeta', t-t'} &  \lrpar{\zeta, \zeta'', t} &  \lrpar{\zeta', \zeta'', t'} } 
\end{equation*}
we have
\begin{equation}\label{comp}
      \curly{K}_{A \circ B} = \lrpar{\beta_{LR,t}}_*
             \lrpar{ \beta_{LM,t-t'}^*\curly{K}_A \cdot \beta_{MR,t'}^* \curly{K}_B } .
\end{equation}

Thus we need a triple heat space, $H^3M_e$, with nice maps ($b$-fibrations) down to the ``left", ``right" and ``center" heat spaces.  In this section, we follow the construction of the heat calculus for closed manifolds in \cite[App. B]{DM} and that of the triple edge space from \cite{Mazzeo:Edge} to construct this triple heat space. For the operators defined above, the resulting formula is given by the following corollary of Theorem \ref{EdgeHeatComp}.

\begin{corollary}\label{CorComp}
Let $A_i \in \Psi^{k_i,\ell_i}_{e, Heat}(M,KD)$ for $i\in \{1,2\}$. 
If $k_i >0$ then the composition is defined, 
      and $A_1 \circ A_2 \in \Psi^{k_1+k_2, \ell_1 + \ell_2}_{e, Heat}(M,KD)$.
\end{corollary}

The construction below works quite generally, as we will discuss at the end of the section. 
As the constructions of this section will not be needed in other sections, the reader may feel free to take the corollary on faith and skip ahead to the next section. 
Throughout this section we consider more general spaces of operators, $\Psi^{\curly{E}}_{e, Heat}$ where $\curly{E}$ is a ``smooth index set" as in section \ref{sec:elliptic}. 
We take as data the existence of a double space and a triple space for the elliptic calculus, along with $b$-fibrations
\begin{equation}\label{Data}
      \xymatrix{ &  M^3_e \ar[dl]^{\beta_{LM}} \ar[d]^{\beta_{LR}} \ar[dr]^{\beta_{MR}} &  \\
         M^2_e &  M^2_e &  M^2_e } .
\end{equation}

It will be useful to start by analyzing the time variables. 
Let $T^2:= \mathbb{R}^+\times \mathbb{R}^+ = \{ (s,s') \}$ and consider the three maps:
\begin{equation*}
      \xymatrix{ &  (s,s') \ar@{|->}[dl]^{\pi_{L}} \ar@{|->}[d]^{\pi_{S}} \ar@{|->}[dr]^{\pi_{R}} &  \\
         s & s+s' & s' }.
\end{equation*}
The first two maps are projections, and easily analyzed. For a function $f$ on $\Rp$ and a density $\mu = g(s,s') ds ds'$ on $T^2$ we have
\begin{equation*}
      \pi_L^*f(s,s') = f(s) , \text{ and } \lrpar{\pi_L}_*\mu = \lrpar{\int g(s,s') ds'} ds.
\end{equation*}
Whereas for the third map,
\begin{equation}\label{pushS} \begin{split}
      \int_{\Rp} f(t) \lrpar{\pi_S}_*\mu &= \int_{T^2} \pi_S^*f(s,s') \mu
      = \int_{T^2} f(s+s') g(s,s') ds ds' \\
      &= \int_{\Rp} f(t) \lrpar{ \int_0^t g(t-t', t') dt' } dt,
\end{split}\end{equation}
hence
\begin{equation*}
      \pi_S^*f(s,s') = f(s+s') , \text{ and } \lrpar{\pi_S}_*\mu = \lrpar{\int_0^s g(s-s',s') ds'} ds.
\end{equation*}
Note that $\pi_S$ is not a $b$-fibration, it is not even a $b$-map. Indeed, $\pi^*(t) = s+s'$ is not a product of boundary defining functions and a positive function. The remedy is to blow-up the corner and consider instead of $T^2$ the space $T^2_0:=[T^2, \{0,0\}]$.
The advantage is that now all of the maps in
\begin{equation*}
      \xymatrix{ & T^2_0:=\lrspar{T^2, \{ \lrpar{0,0} \} }%
             \ar[dl]^{\pi_{L}} \ar[d]^{\pi_{S}} \ar[dr]^{\pi_{R}} &  \\
         \mathbb{R}^+_{t-t'} &  \mathbb{R}^+_{t} &  \mathbb{R}^+_{t'} } 
\end{equation*}
are $b$-fibrations. We will denote the boundary hypersurface at $s=0$ by $\bhs{10}$, that at $s'=0$ by $\bhs{01}$ and the new boundary face from the blow-up, $\bhs{11}$.

Since products of $b$-fibrations are again $b$-fibrations, we set, e.g., $\beta_{MR,R}:=\beta_{MR}\times\pi_{R}$ and obtain
\begin{equation*}
      \xymatrix{ & M^3_e \times T^2_0%
         \ar[dl]^{\beta_{LM,L}} \ar[d]^{\beta_{LR,S}} \ar[dr]^{\beta_{MR,R}}  &  \\
         M^2_e \times \mathbb{R}^+_{t-t'} & M^2_e \times \mathbb{R}^+_{t} &  %
         M^2_e \times \mathbb{R}^+_{t'} } .
\end{equation*}
As we want maps into copies of $HM_e$, we still need to blow-up 
$\lrpar{H\curly{D}_e}_{**, *} := \beta_{**, *}^{-1}\lrpar{ \mathrm{diag}_e \times \{0\} }$ parabolically. 
Similarly, define the ``triple diagonal'', $T\curly{D}_e$, as the closure of the lift of $\{ (m,m,m,0,0) \}$ from $M^3 \times T^2$ to $M^3_e \times T^2_0$.
Define $H^3M_e$ by
\begin{equation*}
      \lrspar{ M^3_e \times T^2_0 ; %
          T\curly{D}_e, \left< ds, ds' \right>;
          \lrpar{H\curly{D}_e}_{LM,L}, \left< ds \right>;
          \lrpar{H\curly{D}_e}_{LR,S}, \left< d(s+s') \right>;
          \lrpar{H\curly{D}_e}_{MR,R}, \left< ds' \right>          }.
\end{equation*}
Recall (\cite{DM} cf. \cite[Prop. 5.12.1]{Corners}) that from a $b$-fibration, $E \xrightarrow{\beta} B$, we get another one by
\begin{equation}\label{CD} \begin{CD}
	 \lrspar{E; \beta^{-1}(S), dt}   @>{\beta_H}>>   E\\
	@V{\tilde{\beta}}VV      @VV{\beta}V \\
	\lrspar{B; S, dt}   @>{\beta_H}>> B
\end{CD} .\end{equation} 

Thus we get $b$-fibrations, e.g., (denoting $t-t'$ by $s$)
\begin{equation*}
\begin{array}{c}
       \lrspar{ M^3_e \times T^2_0 ; %
          \lrpar{T\curly{D}_e}, \left< ds, dt' \right>;
          \lrpar{H\curly{D}_e}_{LM,t-t'}, \left< ds \right>;%
          \lrpar{H\curly{D}_e}_{LR,t}, \left< dt \right>;%
          \lrpar{H\curly{D}_e}_{MR,t'}, \left< dt' \right> } \\
          \downarrow =_{\phantom{LM,t}}\\
                       \lrspar{ M^3_e \times T^2_0 ; %
          \lrpar{H\curly{D}_e}_{LR,t}, \left< dt \right>;%
          \lrpar{T\curly{D}_e}, \left< ds, dt' \right>;
          \lrpar{H\curly{D}_e}_{LM,t-t'}, \left< ds \right>;%
          \lrpar{H\curly{D}_e}_{MR,t'}, \left< dt' \right> } \\ \downarrow \phantom{\beta_{LM,t}}\\
                        \lrspar{ M^3_e \times T^2_0 ; %
          \lrpar{H\curly{D}_e}_{LR,t}, \left< dt \right> }  \\ \downarrow \tilde{\beta}_{LR,t} \\
                        \lrspar{ M^2_e \times \mathbb{R}^+_t ; \mathrm{diag}_e \times \{0\}, \left< dt \right> }
                        =HM_e ,
\end{array} 
\end{equation*}
which fit into the diagram of $b$-fibrations we needed:
\begin{equation}\label{HeatTriple}
      \xymatrix{ & H^3M%
         \ar[dl]^{\beta_{LM,t-t'}} \ar[d]^{\beta_{LR,t}} \ar[dr]^{\beta_{MR,t'}}  &  \\
         \text{  }HM_e \phantom{xx} & \text{  }HM_e\text{  } & \text{  }HM_e\text{  } } .
\end{equation}

The triple heat space has fourteen boundary hypersurfaces. Seven of these, $\bhs{i_1i_2i_3,00}$ with $i_1$, $i_2$, and $i_3$ either $0$ or $1$ (not all $0$), come from the elliptic triple space, see section \ref{sec:elliptic}. Three of these, $\bhs{000,10}$, $\bhs{000,01}$, $\bhs{000,11}$, correspond to $t-t'=0$, $t'=0$, and $t=0$ respectively, away from the diagonals. Another three boundary faces come from blowing up the diagonals: $\bhs{LM,20}$, $\bhs{MR,02}$, and $\bhs{LR,22}$. Finally, there is a boundary face from the triple diagonal, $\bhs{LMR,22}$.

It is useful to have the exponent matrices of the $b$-fibrations in \eqref{HeatTriple}:

{\small
\medskip
\begin{center}
\begin{tabular}{|c|r||c|c|c|c|c|c|c|}      \hline
 & & 001,00 & 010,00 & 100,00 & 110,00 & 101,00 & 011,00 & 111,00  \\  \hline
 & 10,0 & 0 & 0 & 1 & 0 & 1 & 0 & 0 \\ \cline{2-9}
 & 01,0 & 0 & 1 & 0 & 0 & 0 & 1 & 0 \\ \cline{2-9}
 & 11,0 & 0 & 0 & 0 & 1 & 0 & 0 & 1 \\ \cline{2-9}
 & 00,1 & 0 & 0 & 0 & 0 & 0 & 0 & 0 \\ \cline{2-9}
 & 00,2 & 0 & 0 & 0 & 0 & 0 & 0 & 0 \\ \cline{2-9}  \cline{2-9}
 $\beta_{LM,t-t'}$%
 & & 000,01 & 000,10 & 000,11 & MR,02 & LM,20 & LR,22 & LMR,22 \\  \cline{2-9}
 & 10,0 & 0 & 0 & 0 & 0 & 0 & 0 & 0 \\ \cline{2-9}
 & 01,0 & 0 & 0 & 0 & 0 & 0 & 0 & 0 \\ \cline{2-9}
 & 11,0 & 0 & 0 & 0 & 0 & 0 & 0 & 0 \\ \cline{2-9}
 & 00,1 & 0 & 1 & 1 & 0 & 0 & 2 & 0 \\ \cline{2-9}
 & 00,2 & 0 & 0 & 0 & 0 & 1 & 0 & 1 \\ \hline
\end{tabular}
\end{center} }

{\small
\medskip
\begin{center}
\begin{tabular}{|c|r||c|c|c|c|c|c|c|}      \hline
 & & 001,00 & 010,00 & 100,00 & 110,00 & 101,00 & 011,00 & 111,00  \\  \hline
 & 10,0 & 0 & 0 & 1 & 1 & 0 & 0 & 0 \\ \cline{2-9}
 & 01,0 & 1 & 0 & 0 & 0 & 0 & 1 & 0 \\ \cline{2-9}
 & 11,0 & 0 & 0 & 0 & 0 & 1 & 0 & 1 \\ \cline{2-9}
 & 00,1 & 0 & 0 & 0 & 0 & 0 & 0 & 0 \\ \cline{2-9}
 & 00,2 & 0 & 0 & 0 & 0 & 0 & 0 & 0 \\ \cline{2-9}  \cline{2-9}
 $\beta_{LR,t}$%
 & & 000,01 & 000,10 & 000,11 & MR,02 & LM,20 & LR,22 & LMR,22 \\  \cline{2-9}
 & 10,0 & 0 & 0 & 0 & 0 & 0 & 0 & 0 \\ \cline{2-9}
 & 01,0 & 0 & 0 & 0 & 0 & 0 & 0 & 0 \\ \cline{2-9}
 & 11,0 & 0 & 0 & 0 & 0 & 0 & 0 & 0 \\ \cline{2-9}
 & 00,1 & 0 & 0 & 1 & 0 & 0 & 0 & 0 \\ \cline{2-9}
 & 00,2 & 0 & 0 & 0 & 0 & 0 & 1 & 1 \\ \hline
\end{tabular}
\end{center}}

{\small
\medskip
\begin{center}
\begin{tabular}{|c|r||c|c|c|c|c|c|c|}      \hline
 & & 001,00 & 010,00 & 100,00 & 110,00 & 101,00 & 011,00 & 111,00  \\  \hline
 & 10,0 & 0 & 1 & 0 & 1 & 0 & 0 & 0 \\ \cline{2-9}
 & 01,0 & 1 & 0 & 0 & 0 & 1 & 0 & 0 \\ \cline{2-9}
 & 11,0 & 0 & 0 & 0 & 0 & 0 & 1 & 1 \\ \cline{2-9}
 & 00,1 & 0 & 0 & 0 & 0 & 0 & 0 & 0 \\ \cline{2-9}
 & 00,2 & 0 & 0 & 0 & 0 & 0 & 0 & 0 \\ \cline{2-9}  \cline{2-9}
 $\beta_{MR,t'}$%
 & & 000,01 & 000,10 & 000,11 & MR,02 & LM,20 & LR,22 & LMR,22 \\  \cline{2-9}
 & 10,0 & 0 & 0 & 0 & 0 & 0 & 0 & 0 \\ \cline{2-9}
 & 01,0 & 0 & 0 & 0 & 0 & 0 & 0 & 0 \\ \cline{2-9}
 & 11,0 & 0 & 0 & 0 & 0 & 0 & 0 & 0 \\ \cline{2-9}
 & 00,1 & 1 & 0 & 1 & 0 & 0 & 2 & 0 \\ \cline{2-9}
 & 00,2 & 0 & 0 & 0 & 1 & 0 & 0 & 1 \\ \hline
\end{tabular}
\end{center} }

Next, we will rewrite \eqref{comp} in terms of $b$-densities, as we can then apply a convenient form of the push-forward theorem.
We denote by $\nu_X$ a standard non-vanishing half-density on $X$, and abbreviate $\nu_{HM_e}$ and $\nu_{H^3M}$ to $\nu_H$ and $\nu_{H^3}$ respectively.
Note that
$\curly{K}_A = \kappa_A \rho_{00,2}^{-{\frac{m+3}{2}}} \rho_{11,0}^{-{\frac{b+1}{2}}} \nu_H$,
so \eqref{comp} becomes
\begin{equation*}
            \kappa_{A \circ B} \rho_{00,2}^{-{\frac{m+3}{2}}} \rho_{11,0}^{-{\frac{b+1}{2}}} \nu_H =
             \lrpar{\beta_{LR,S}}_*
             \lrspar{ \beta_{LM,L}^*\lrpar{ %
                \kappa_A \rho_{00,2}^{-{\frac{m+3}{2}}} \rho_{11,0}^{-{\frac{b+1}{2}}} \nu_H } 
              \cdot \beta_{MR,R}^*\lrpar{ %
                 \kappa_B \rho_{00,2}^{-{\frac{m+3}{2}}} \rho_{11,0}^{-{\frac{b+1}{2}}} \nu_H } } .
\end{equation*}
Abbreviating pull-back via subindices and multiplying both sides by
$\rho_{00,2}^{-{\frac{m+3}{2}}} \rho_{11,0}^{-{\frac{b+1}{2}}} \nu_H$, this becomes
\begin{equation}\label{comp2}
            \kappa_{A \circ B} \rho_{00,2}^{-\lrpar{m+3}} \rho_{11,0}^{-\lrpar{b+1}} \nu_H^2 =
             \lrpar{\beta_{LR,S}}_*
             \lrspar{ \lrpar{\kappa_A}_{LR} \lrpar{\kappa_B}_{MR} 
             \lrpar{\rho_{00,2}^{-{\frac{m+3}{2}}} \rho_{11,0}^{-{\frac{b+1}{2}}} \nu_H }_{LR,LM,MR} }.
\end{equation}

This is easily computed making use of the commutative diagram(s):
\begin{equation}\label{BigDig}
      \xymatrix{ H^3M \ar@/^/[ddd] \ar[ddd] \ar@/_/[ddd] \ar[rr]^{\beta_{H^3}} \ar[rd]_{\gamma_1} &  & %
         M^3 \times T^2_0 \ar@/^/[ddd] \ar[ddd] \ar@/_/[ddd] \\
         & M^3_e \times T^2_0 \ar[ru]_{\gamma_2} & \\
         & M^2_e \times \mathbb{R}^+ \ar[rd] & \\
         HM_e \ar[ru]  \ar[rr]^{\beta_H} &  & M^2 \times \mathbb{R}^+ },
\end{equation}
where the left-most column is the diagram \eqref{HeatTriple}.
Indeed, a consequence of \eqref{pushS} is that 
\begin{equation*}
      \lrpar{\pi_S}_*(ds ds') = t dt
\end{equation*}
or $\pi_S^*(\Omega_{\Rp}) = \rho_{11}^{-1}\Omega_{T^2_0}$.
Using this, it is easily seen that
\begin{equation*}
      \lrpar{\nu_{M^2\times \Rp}}_{LR,LM,MR} = \nu_{M^3 \times T^2_0}^2 ,
\end{equation*}
and that (cf. \cite[Lemma 2.2]{Vaillant})
{\small
\begin{equation}\label{DensityPullBack}\begin{split}
      \beta_{H}^*&\lrpar{ \Omega_{M^2 \times \Rp} } =
         \rho_{00,2}^{m+1} \rho_{11,0}^{b+1}  \Omega_{HM}, \\
      \beta_{H^3}^*&\lrpar{ \Omega_{M^3 \times T^2_0} } =
         \gamma_1^*\lrpar{ \lrpar{ \rho_{110,00}\rho_{101,00} \rho_{011,00} }^{b+1}%
            \rho_{111,00}^{2b+2} \Omega_{M_e^3 \times T^2_0} } \\
         &=  \lrpar{ \rho_{110,00}\rho_{101,00} \rho_{011,00} }^{b+1}%
            \rho_{111,00}^{2b+2}
            \gamma_1^*\lrpar{ \Omega_{M_e^3 \times T^2_0} } \\
         &=  \lrpar{ \rho_{110,00}\rho_{101,00} \rho_{011,00} }^{b+1}%
            \rho_{111,00}^{2b+2}
            \lrpar{ \rho_{MR,02} \rho_{LM,20} \rho_{LR,22} }^{m+1} \rho_{LMR,22}^{2m+3}
            \Omega_{H^3M},
\end{split}\end{equation} }
so that
{\small
\begin{equation*}
       \begin{split}
       \left(\rho_{00,2}^{-{\frac{m+3}{2}}} \right.&\left.\rho_{11,0}^{-{\frac{b+1}{2}}} \nu_H \right)_{LR,LM,MR} 
        = \lrspar{ \rho_{00,2}^{-{\frac{m+3}{2}}} \rho_{11,0}^{-{\frac{b+1}{2}}} %
              \lrpar{ \rho_{00,2}^{-\frac{m+1}{2}} \rho_{11,0}^{-\frac{b+1}{2}} %
              \beta_H^*\lrpar{ \nu_{M^2 \times \Rp}} } }_{LM,MR,LR} \\
        &= \lrpar{ \rho_{00,2}^{-\lrpar{m+2}} \rho_{11,0}^{-\lrpar{b+1}} }_{LM,MR,LR} %
            \beta_{H^3}^*\lrpar{ \nu_{M^3 \times T^2_0}^2} \\
        &= \lrpar{ \rho_{MR,02} \rho_{LM,20} \rho_{LR,22} \rho_{LMR,22}^3 }^{-\lrpar{m+2}} 
              \lrpar{ \rho_{110,00}\rho_{101,00} \rho_{011,00} \rho_{111,00}^3 }^{-\lrpar{b+1}} \\
        &\phantom{x}
              \left[ \lrpar{ \rho_{110,00}\rho_{101,00} \rho_{011,00} }^{b+1}%
              \rho_{111,00}^{2b+2} %
              \lrpar{ \rho_{MR,02} \rho_{LM,20} \rho_{LR,22} }^{m+1} 
              \rho_{LMR,22}^{2m+3}  \nu_{H^3M}^2 \right] \\
        &= \lrpar{ \rho_{MR,02} \rho_{LM,20} \rho_{LR,22} }^{-1} \rho_{LMR,22}^{-m-3}
              \rho_{111,00}^{-b-1} \nu_{H^3M}^2 . \\
\end{split}  \end{equation*} }
Notice that thus \eqref{comp2} is
{\small
\begin{equation}\label{comp3} \begin{split}
            \kappa_{A \circ B} &\rho_{00,2}^{-m-3}  \rho_{11,0}^{-b-1} \nu_H^2  
            = \lrpar{\beta_{LR,t}}_* \lrspar{ \lrpar{\kappa_A}_{LM} \lrpar{\kappa_B}_{MR} %
              \lrpar{ \rho_{MR,02} \rho_{LM,20} \rho_{LR,22} }^{-1} \rho_{LMR,22}^{-m-3}%
              \rho_{111,00}^{-b-1} \nu_{H^3M}^2 } .
\end{split} \end{equation} }
We change this to $b$-densities. Recall that a $b$-density is a smooth density divided by a ``total" boundary defining function, hence
\begin{gather*}
      \nu^2_H = \lrpar{\rho_{10,0}\rho_{01,0}\rho_{11,0}\rho_{00,1}\rho_{00,2}} {}^b\nu_H^2 %
      \text{ and }\\
      \nu^2_{H^3} = \lrpar{\rho_{001,00} \cdots
            %\rho_{010,00} \rho_{100,00} \rho_{110,00} \rho_{101,00} \rho_{011,00} 
            \rho_{111,00} \rho_{000,01} \rho_{000,10} \rho_{000,11}%
            \rho_{MR,02} \rho_{LM,20} \rho_{LR,22} \rho_{LMR,22}} {}^b\nu_{H^3}^2 .
\end{gather*}
Substituting into \eqref{comp3} and simplifying yields
{\small
\begin{equation}\label{comp4}
            \kappa_{A \circ B} {}^b\nu_H^2 
            = \lrpar{\beta_{LR,t}}_* \lrspar{ \lrpar{\kappa_A}_{LM} \lrpar{\kappa_B}_{MR} %
               \rho_{101,00}^{b+1} \rho_{LR,22}^{m+2} \rho_{010,00} \rho_{000,01} 
               \rho_{000,10} {}^b\nu_{H^3M}^2 }.
\end{equation} 
Finally, we can conclude using the push-forward theorem.

\begin{theorem}\label{EdgeHeatComp}
Let $A \in \Psi^\curly{E}_{e, Heat}$ and  $B \in \Psi^\curly{F}_{e, Heat}$. Then provided 
\begin{itemize}
       \item $\Re(E_{01,0})+\Re(F_{10,0})+1>0$
       \item $\Re(E_{00,1})+1>0$
       \item $\Re(F_{00,1})+1>0$
       \item $\Re(E_{00,2})>0$
       \item $\Re(F_{00,2})>0$,
\end{itemize}
the composition $A \circ B$ is well-defined and is an element of $\Psi^\curly{G}_{e, Heat}$ where
\begin{itemize}
       \item $G_{10,0} = E_{10,0} \bar{\cup} \lrpar{E_{11,0} + F_{10,0}}$
       \item $G_{01,0} = F_{01,0} \bar{\cup} \lrpar{E_{01,0} + F_{11,0}}$
       \item $G_{11,0} = \lrpar{E_{11,0} + F_{11,0}} \bar{\cup} \lrpar{E_{10,0} + F_{01,0}+b+1}$
       \item $G_{00,1} = E_{00,1} + F_{00,1} $
       \item $G_{00,2} = \lrpar{E_{00,2} + F_{00,2}} \bar{\cup} \lrpar{2E_{00,1} + 2F_{00,1} + m + 2}$.
\end{itemize}
\end{theorem}

\begin{proof}
We apply the push-forward theorem using the exponent matrices and \eqref{comp4}. The integrability conditions are at the faces: $\bhs{010,00}$, $\bhs{000,01}$, $\bhs{000,10}$, $\bhs{000,02}$, and $\bhs{000,20}$, since these are mapped into the interior by $\beta_{LR,S}$.
\end{proof}

Looking back over the proof of the theorem, it is easy to see that the same construction works much more generally. The point is that for all Metrics with Iterated Complete Edge structures, or MICE,
the heat space is constructed in the same manner as $HM_e$. Namely take the appropriate double space, say $M_{\text{any}}^2$, and let 
\begin{equation*}
      HM_{\text{any}} 
         = [M_{\text{any}}^2\times \mathbb{R}^+; \text{diag}_{\text{any}} \times \{0\}, \left< dt \right> ]. 
\end{equation*}
Similarly, the above construction of the heat triple space plays out the same way given only the diagram \eqref{Data}.
As for the densities, note that the construction of the triple heat space factors through 
$M^3_{\text{any}} \times T^2_0$ as in \eqref{BigDig} with all of the remaining blow-ups involving only the temporal variables.

So assume that we have such a calculus, $\Psi^{\curly{E}}_{\text{any}, Heat}$. Denote by $\curly{E}'$ the part of $\curly{E}$ corresponding to spatial boundary hypersurfaces and $\curly{E}'' = (E''_1, E''_2)$ the index sets corresponding to the temporal boundaries.

\begin{theorem}
Let $A \in \Psi^\curly{E}_{\text{any}, Heat}$ and  $B \in \Psi^\curly{F}_{\text{any}, Heat}$. Assume that operators in $\Psi^\curly{E'}_{\text{any}}$ and $\Psi^\curly{F'}_{\text{any}}$ compose with resulting index set
\begin{equation*}
      \Psi^\curly{E'}_{\text{any}} \circ \Psi^\curly{F'}_{\text{any}} \subseteq \Psi^{\curly{G'}}_{\text{any}}.
\end{equation*}
Then provided 
\begin{itemize}
       \item $\Re(E''_{1})+1>0$
       \item $\Re(F''_{1})+1>0$
       \item $\Re(E''_{2})>0$
       \item $\Re(F''_{2})>0$,
\end{itemize}
the composition $A \circ B$ is well-defined and is an element of $\Psi^\curly{G}_{\text{any}, Heat}$ where $\curly{G}'$ is as above and 
\begin{itemize}
       \item $G''_{1} = E''_{1} + F''_{1} $
       \item $G''_{2} = \lrpar{E''_{2} + F''_{2}} \bar{\cup} \lrpar{2E''_{1} + 2F''_{1} + m + 2}$.
\end{itemize}
\end{theorem}

That is to say, for MICE, the heat calculus is just the elliptic calculus together with the heat calculus for closed manifolds. Notice that the situation is different for {\em incomplete} metrics. Heuristically, the ``heat" arrives at the boundary in finite time, and thus the boundary (or parts of it) need to be blown-up at $\{ t=0 \}$. For instance the heat space for an incomplete edge metric is pictured below.

\begin{figure}[htpb]
      \centering
      \includegraphics[bb=0in 2.25in 8.5in 8.75in, keepaspectratio, height=2in,clip]{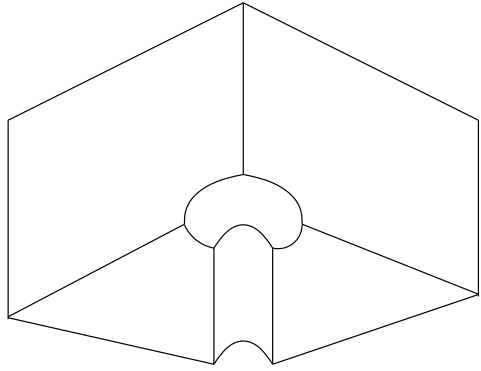}
      \caption{Heat space for an incomplete metric}
      \label{DHeatSpace}
\end{figure}

%%%%%%%%%%%%%%%%%%%%%%%%%%%%%%%%%%%%%%%%%%%%%%%%%%%%%%%%%%%%%%%%%%%%%%%%%%%%%%%%%%%%%%%%%%%%%%%%%%%%%%%%%%%%%%%%%%%%%%%%%%%%
\section{Construction of the Heat Kernel} \label{sec:Construct}
%%%%%%%%%%%%%%%%%%%%%%%%%%%%%%%%%%%%%%%%%%%%%%%%%%%%%%%%%%%%%%%%%%%%%%%%%%%%%%%%%%%%%%%%%%%%%%%%%%%%%%%%%%%%%%%%%%%%%%%%%%%%

In this section we will construct a solution to the heat equation 
by constructing its Schwartz kernel as a distribution on $HM_e$. The plan is to solve away its Taylor series at each of the boundary hypersurfaces and then deal with the remainder. To this end, we shall make persistent use of the ``normal operators".
These are extra symbol maps necessitated by the presence of the boundary. They are given by restricting the suitably weighted kernel of the operator to each of the boundary hypersurfaces. Thus if $A \in \Psi^{k, \ell}_{e, Heat}(M,KD)$, then its normal operators are
\begin{equation}
      N_{11,0}^{\ell}(A) = \rho_{11,0}^{-\ell} \curly{K}_A \downharpoonright_{\bhs{11,0}}
      \text{ and }
      N_{00,2}^{k}(A) = \rho_{00,2}^{-k} \curly{K}_A \downharpoonright_{\bhs{00,2}}.
\end{equation}
Just like the symbol map in the usual pseudodifferential calculus, these maps 
fit into short exact sequences:
\begin{equation}\label{SES} \begin{split}
      &0 \to \Psi^{k-1, \ell}_{e, Heat} \to \Psi^{k, \ell}_{e, Heat} \xrightarrow{N_{00,2}} 
           \rho^{\ell}_{11,0}  \dot{\curly{C}}^{\infty}_{11,0}
            \lrpar{\bhs{00,2}, KD\downharpoonright_{\bhs{00,2}}} \to 0 \\
      &0 \to \Psi^{k, \ell-1}_{e, Heat} \to \Psi^{k, \ell}_{e,Heat} \xrightarrow{N_{11,0}} 
            \rho^{k}_{00,2}  \dot{\curly{C}}^{\infty}_{00,2}
            \lrpar{\bhs{11,0}, KD\downharpoonright_{\bhs{11,0}}} \to 0.
\end{split} \end{equation}
The spaces on the right can in turn be interpreted as spaces of kernels of operators from simpler calculi. Indeed, $\bhs{00,2}$ fibers over the diagonal of $M$ and can be thought of as a compactified (edge) tangent bundle and the normal operator at $\bhs{11,0}$ can be thought of as a family of normal operators from the elliptical calculus.

The heat equation, when restricted to the boundary hypersurfaces, induces equations for the corresponding normal operators. More generally, we have the following lemma wherein we will think of the symbol of an (edge) differential operator as a constant coefficient differential operator on the fibers of the (edge) tangent bundle.
\begin{lemma}
Let $A \in \Psi^{k, \ell}_{e, Heat}(M, KD)$. If $V$ is any edge vector field on $M$, with symbol ${}^e\sigma(V)$ and normal operator $N_e(V)$, then $t^{1/2}V \circ A \in \Psi^{k,\ell}_{e, Heat}(M,KD)$ and 
\begin{equation*}
      N_{11,0}^{\ell} \lrpar{t^{1/2}V \circ A} = t^{1/2}N_e(V)N_{11,0}^{\ell}(A), \text{ and }
      N_{00,2}^{k} \lrpar{t^{1/2}V \circ A} = {}^e\sigma(V)N_{00,2}^{k}(A).      
\end{equation*}
Similarly, $t\partial_t \circ A \in \Psi^{k,\ell}_{e, Heat}(M,KD)$ with normal operators
\begin{equation*}
      N_{11,0}^{\ell} \lrpar{t\partial_t \circ A} = t\partial_tN_{11,0}^{\ell}(A), \text{ and }
      N_{00,2}^{k} \lrpar{t\partial_t \circ A} = -\frac{1}{2}\lrpar{\curly{R}+m-k+2}N_{00,2}^{k}(A).      
\end{equation*}
Here $\curly{R}$ denotes the radial vector field on the fibers of ${}^eTM$.
\end{lemma}
\begin{proof}
The proof consists of a computation in local coordinates. We will use coordinates as in \eqref{ConvCoord} 
\begin{equation}
      \lrpar{ \zeta, \lrpar{ \frac{x-x'}{xt^{1/2}}, \frac{y-y'}{xt^{1/2}}, \frac{z-z'}{t^{1/2}} }, t^{1/2} }=:
          \lrpar{ \zeta, \lrpar{ \curly{S}', \curly{U}', \curly{Z}'}, \tau} .
\end{equation}
Let $V$ be an edge vector field, and $V'$ its adjoint. Then we have
\begin{equation}\label{KerVA}
      \left< V\widetilde{A}\psi , \phi \right> 
      = \left< \widetilde{A}\psi, V'\phi \right>
      = \left< \curly{K}_A, \beta_H^*\lrpar{V '\phi \hat{*}_t \psi} \right>
      = \left< \beta_{H,L}^*(V')'  \curly{K}_A, \beta_H^*\lrpar{\phi \hat{*}_t \psi } \right>,
\end{equation}
where $\beta_{H,L}$ is the map $HM_e \to M^2_e \times \Rp \to M$ induced by projection onto the left factor of $M$. Notice that the integration by parts implicit in \eqref{KerVA} are justified since 
%there is no derivative with respect to $t$ and 
edge vector fields are tangent to all spatial boundary faces.
It is important to keep track of the densities in \eqref{KerVA}, note that by \eqref{DenisityPullBack} and \eqref{Jacobian} for any $\Psi \in C^{\infty}\lrpar{M \times \mathbb{R}^+;\Omega^{1/2}}$,
\begin{equation*}\begin{split}
	\left<\curly{K}_A %
	\left\lvert \frac{d\zeta d\curly{S'} d\curly{U'} d\curly{Z'} d\tau}%
	{\tau^{\frac{m+3}2-k}x^{\frac{b+1}2-\ell}} \right\rvert^{\frac12},%
	\beta_H^*\lrpar{ \Psi \left\lvert d\zeta d\zeta' dt \right\rvert^{\frac12}} \right>
	&=
	\left<\curly{K}_A %
	\left\lvert \frac{d\zeta d\curly{S'} d\curly{U'} d\curly{Z'} d\tau}%
	{\tau^{\frac{m+3}2-k}x^{\frac{b+1}2-\ell}} \right\rvert^{\frac12},%
	\beta_H^*\lrpar{ \Psi } %
	\left\lvert \frac{d\zeta d\curly{S'} d\curly{U'} d\curly{Z'} d\tau}%
	{\tau^{-\frac{m+1}2}x^{-\frac{b+1}2}} \right\rvert^{\frac12} \right>\\%
	&=
	\left<\curly{K}_A %
	\tau^{k-1}x^{\ell}%
	\left\lvert d\zeta d\curly{S'} d\curly{U'} d\curly{Z'} d\tau \right\rvert^{\frac12}, %
	\beta_H^*\lrpar{ \Psi } %
	\left\lvert d\zeta d\curly{S'} d\curly{U'} d\curly{Z'} d\tau \right\rvert^{\frac12} \right> \\
\end{split}\end{equation*}
Thus the kernel of $VA$ is $\beta_{H,L}^*(V')' \curly{K}_A$, and we can find the normal operators by
conjugating with a defining function.
For instance for $W$ a vector field as above, at $\bhs{11,0}$,
\begin{equation*}
	N_{11,0}\lrpar{ W \circ A} = x^{-\ell} \curly{K}_{W\circ A}\downharpoonright_{11,0}
	=N_{11,0}(W) \lrpar{x^{-\ell} \curly{K}_{A}}\downharpoonright_{11,0},
\end{equation*}
so $N_{11,0}(W) = x^{-\ell} \beta_{H,L}^*(W')' x^{\ell}\downharpoonright_{11,0}$,
while similar reasoning yields
$N_{00,2}(W) = \tau^{-(k-1)} \beta_{H,L}^*(W')' \tau^{k-1}\downharpoonright_{00,2}$
(the $k-1$ from the densities).

We can carry out the lifts:
\begin{gather*}
      \beta_{H,L}^*\lrpar{x\partial_x} = x\lrspar{%
           \partial_x + \frac{x'}{x^2t^{1/2}} \partial_{\curly{S}'} - \frac{y-y'}{x^2t^{1/2}} \partial_{\curly{U}'}}
           = x\partial_x - \curly{S}' \partial_{\curly{S}'} - \curly{U}' \partial_{\curly{U}'} + 
           \frac{\partial_{\curly{S}'}}{t^{1/2}} \\
      \beta_{H,L}^*\lrpar{x\partial_y} = x \lrspar{%
           \partial_y + \frac{1}{xt^{1/2}}\partial_{\curly{U}'} } = x\partial_y + \frac{\partial_{\curly{U}'}}{t^{1/2}}\\
      \beta_{H,L}^*\lrpar{\partial_z} = \partial_z + \frac{\partial_{\curly{Z}'}}{t^{1/2}}.
\end{gather*}
Hence $N_{00,2}(t^{1/2}VA) = {}^e\sigma\lrpar{V}N_{00,2}(A)$ as we only keep that part of the lift divided by $t^{1/2}$. Note that $N_{11,0}(t^{1/2}VA) = t^{1/2}N(V)N_{11,0}(A)$ with $N(V)$ the normal operator from the elliptic calculus as near $\bhs{11,0}$ and away from $\bhs{00,2}$ we can use coordinates in which the temporal and spatial variables do not interact, e.g. \eqref{efCoord}.

Similarly, consider the effect of $t\partial_t$. Its lift under the projection $HM_e \to M^2_e \times \Rp \to \Rp$ is given by
\begin{equation*}
	\tau^2\lrspar{ \frac{1}{2\tau}\partial_{\tau} - \frac{1}{2\tau^2}\lrpar{\curly{S}'\partial_{\curly{S}'} +%
	\curly{U}'\partial_{\curly{U}'} + \curly{Z}'\partial_{\curly{Z}'} } }
	= \frac{\tau}{2}\partial_{\tau} - \frac{1}{2} \curly{R} ,
\end{equation*}
so since
\begin{equation*}\begin{split}
	\tau^{-k+1} \beta_{H,L}^*&\lrpar{ \lrpar{t\partial_t}'}' \tau^{k-1}
	=\tau^{-k+1} \beta_{H,L}^*(-Id - t\partial_t)' \tau^{k-1}
	=-Id -  \tau^{-k+1} \lrpar{ \frac{\tau}{2}\partial_{\tau} - \frac{1}{2} \curly{R} }' \tau^{k-1} \\
	&=-Id +  \tau^{-k+1}
		\lrpar{ \frac{Id}{2} + \frac{\tau}{2}\partial_{\tau} - \frac{1}{2} \lrpar{ m +\curly{R} } } \tau^{k-1}
%	=-Id + \frac{Id}{2} + \frac{k-1}{2}Id + \frac{\tau}{2}\partial_{\tau}- \frac{1}{2} \lrpar{ m +\curly{R} }
	=\frac{\tau}{2}\partial_{\tau}- \frac{1}{2} \lrpar{\curly{R} + m-k+2 }
\end{split}\end{equation*}
we conclude that $N_{00,2}^k(t\partial_tA) = - \frac{1}{2} \lrpar{ \curly{R} + m - k + 2}N_{00,2}^k(A)$. 
\end{proof}

This lemma, together with the composition result from section $\mathcal{x}$\ref{sec:Comp}, allows us to construct a solution to the heat equation.

\begin{theorem}\label{HeatKerExt}
Let $\Delta$ be the Laplacian of an exact edge metric.
There exists a unique solution $\curly{H} \in \Psi^{2,0}_{e,Heat}(M,KD)$ to the heat equation
\begin{equation}\label{HeatEq}
      \begin{cases}
          \lrpar{\partial_t + \Delta}\curly{H} = 0\\
          \displaystyle \lim_{t \to 0}\curly{H} = \text{Id}
      \end{cases}
      \iff  \lrpar{\partial_t + \Delta} \widetilde{\curly{H}} = \widetilde{\text{Id}} .
\end{equation}
Moreover, $\curly{H}$ has normal operators:
\begin{gather}
      \label{NormOp1}
            N^{2}_{00,2}\lrpar{\curly{H}} = \frac{1}{\lrpar{4\pi}^{\frac{n}{2}} }
                  \exp\lrpar{-\frac{1}{4}|v|^2_e}\text{Id} \text{ at $(\zeta,v) \in {}^eT^*M$} \\
      \label{NormOp2}
            N^{0}_{11,0}\lrpar{\curly{H}} = \exp\lrpar{-t\Delta_{\mathbb{H}^b}} \exp\lrpar{-t\Delta_{F}}
\end{gather}
\end{theorem}

\begin{proof}
By the previous lemma, the heat equation imposes
\begin{gather}
     \label{NormalEq}
            \lrpar{ {}^e\sigma(\Delta) - \frac{1}{2}\lrpar{\curly{R}+m} }N_{00,2}^2(\curly{H}) = 0 \text{ and } \\ 
     \label{NormalEq2}
            \lrpar{ \partial_t + N_e\lrpar{\Delta} } N^0_{11,0} \lrpar{ \curly{H} } = 0.
\end{gather}
Furthermore, we have a boundary condition, most easily expressed using \eqref{EffectiveAct}. Indeed, note that if $A \in \Psi^{2,0}_{e,Heat}$ then, upon restricting to $t=0$, $A$ acts as a multiplication operator
\begin{equation*}
      Af(\zeta)= \lrspar{ \int%
              \kappa_A\lrpar{\zeta, \curly{S}',\curly{U}',\curly{Z}', 0}%
              d\curly{S}' d\curly{U}' d\curly{Z}' }
              f\lrpar{\zeta} |d\zeta dt|^{1/2} .
\end{equation*}
Our boundary condition is thus
\begin{equation}\label{BdyCon}
              \int \kappa_A\lrpar{\zeta, \curly{S}',\curly{U}',\curly{Z}', 0}%
              d\curly{S}' d\curly{U}' d\curly{Z}' =1 \iff
              \int_{\text{fiber}} N^2_{00,2}\lrpar{\curly{H}} = \text{Id}.      
\end{equation}
Equations \eqref{NormalEq} and \eqref{BdyCon} are fibre-by-fibre conditions. As in \cite{APS Book}, on any fixed fiber we can choose coordinates so that ${}^e\sigma(\Delta)$ is the Laplacian on $\mathbb{R}^n$, the unique solution is then seen to be \eqref{NormOp1}.

Clearly, \eqref{NormOp2} solves \eqref{NormalEq2} and \eqref{NormOp1}, \eqref{NormOp2} are consistent. Hence there exists $G_{(1)} \in \Psi^{2,0}_{e, Heat}(M,KD)$ with normal operators given by \eqref{NormOp1}, \eqref{NormOp2} and satisfying \eqref{NormalEq}, \eqref{NormalEq2}, \eqref{BdyCon}. By exactness of the sequences \eqref{SES}, this implies that
\begin{equation*}
      \lrpar{\partial_t + \Delta} \widetilde{G_{(1)}} = \widetilde{\text{Id}} - \widetilde{R_{(1)}},
\end{equation*}
with  $R_{(1)} \in \Psi^{3,1}_{e, Heat}(M,KD)$.
By the composition formula, Corollary \ref{CorComp}, $R^k \in \Psi^{3k, k}_{e, Heat}$ hence the series $\sum R^k$ can be summed, say to $\text{Id} + S$, $S \in \Psi^{3,1}_{e, Heat}$.
Hence, with 
\begin{equation*}
      \widetilde{G_{(2)}}= \widetilde{G_{(1)}} \circ \lrpar{\widetilde{\text{Id}} + \widetilde{S}},
\end{equation*}
\begin{equation*}
      \lrpar{\partial_t + \Delta} \widetilde{G_{(2)}} = \widetilde{\text{Id}} - \widetilde{R_{(2)}},
\end{equation*}
with  $R_{(2)} \in \Psi^{\infty,\infty}_{e, Heat}(M,KD)$. Thus $G_{(1)}$ and $G_{(2)}$ are parametrices of first order and infinite order, respectively. Finally, any element of $\text{Id} + \Psi^{\infty, \infty}_{e, Heat}$ has an inverse in the same space (cf. \cite[Prop. 7.17]{APS Book}), so $G_{(2)}$ 
can be improved to an actual inverse, $\curly{H}$.
\end{proof}

Just as in \cite[$\mathcal{x}$7.6]{APS Book}, one can extend the discussion to generalized Dirac operators with bundle coefficients. If $E$ is a bundle over $M$, we define the space of operators just as before but with the coefficient bundle replaced by 
\begin{equation}\label{BundleCoeff}
      \mathscr{K}^{k,\ell}(M,E) =  \mathscr{K}^{k,\ell}(M,KD) \otimes_{C^{\infty}(HM_e)} 
            C^{\infty}\lrpar{HM_e; \beta_H^*\lrpar{\Hom\lrpar{E \otimes KD^*}}}.
\end{equation}
The discussion of composition and the construction of the heat kernel can be extended to this context.

%%%%%%%%%%%%%%%%%%%%%%%%%%%%%%%%%%%%%%%%%%%%%%%%%%%%%%%%%%%%%%%%%%%%%%%%%%%%%%%%%%%%%%%%%%%%%%%%%%%%%%%%%%%%%%%%%%%%%%%%%%%%
\section{The Index Theorem} \label{sec:IndexTheorem}
%%%%%%%%%%%%%%%%%%%%%%%%%%%%%%%%%%%%%%%%%%%%%%%%%%%%%%%%%%%%%%%%%%%%%%%%%%%%%%%%%%%%%%%%%%%%%%%%%%%%%%%%%%%%%%%%%%%%%%%%%%%%

On a closed manifold McKean and Singer noticed the remarkable fact that the super-trace of the heat kernel is independent of $t$. Since the index is the limit at infinity of the heat kernel, one obtains the index theorem from the short-time asymptotics of the heat kernel. On an edge manifold, the heat kernel is unfortunately not trace-class. Nevertheless, we obtain an index theorem by renormalizing the super-trace of the heat kernel and comparing its values as time goes to zero and infinity.

%%%%%%%%%%%%%%%%%%%%%%%%%%%%%%%%%%%%%%%%%%%%%%%%%%%%%%%%%%%%%%%%%%%%%%%%%%%%%%%%%%%%%%%%%%%%%%%%%%%%%%
\subsection{Renormalization} \label{sec:RInt} 
%%%%%%%%%%%%%%%%%%%%%%%%%%%%%%%%%%%%%%%%%%%%%%%%%%%%%%%%%%%%%%%%%%%%%%%%%%%%%%%%%%%%%%%%%%%%%%%%%%%%%%

%%%%%%%%%%%%%%%%%%%%%%%%%%%%%%%%%%%%%%%%%%%%%%%%%%%%%%%%%%%%%%%%%%%%%%%%%%%%%%%%%%%%%%%%%%%%%%%%%%%%%%
\subsubsection{Renormalized Integrals} \label{subsec:RInt} $ $\newline
%%%%%%%%%%%%%%%%%%%%%%%%%%%%%%%%%%%%%%%%%%%%%%%%%%%%%%%%%%%%%%%%%%%%%%%%%%%%%%%%%%%%%%%%%%%%%%%%%%%%%%

A manifold with an edge metric, $(M,g)$ is topologically a manifold with boundary, and the study of edge metrics naturally involves densities defined on $M$ with expansions at the boundary, of the form
\begin{equation*}
      \mu \sim \lrspar{\sum_{k>-N} \sum_{j=0}^{n_k} a_{k,j}(y,z) x^k \lrpar{ \log x }^j } dx \dvol_{\partial M},
\end{equation*}
with $a_{k,j}$ smooth. Expansions of this form are known as polyhomogeneous conormal distributions, and will be referred to as phg. Naturally, the coefficients in this expansion depend on the choice of bdf, but the existence of the expansion and some information about the exponents involved is independent of this choice. Densities of this form are generally not integrable, yet it is precisely their integral that will interest us!

Recall how the $\Gamma$ function is defined on the complex plane.
We start with an explicit integral expression
\begin{equation*}
      \Gamma(z) = \int_0^{\infty} t^{z-1}e^{-t} dt
\end{equation*}
defined on a half plane $\{\mathrm{Re}(z) > 0 \}$. Then we use the series expansion of the exponential to meromorphically continue $\Gamma(z)$ to the plane.

Similarly, given a density $\mu$ with a phg expansion and a choice of bdf, $x$ , we consider the zeta function
\begin{equation*}
      \zeta_x(z) = \int x^z \mu.
\end{equation*}
This is initially defined on a half plane, but can be meromorphically extended to the complex plane by using the phg expansion of $\mu$. We define the renormalized integral of $\mu$ as the finite part of $\zeta_x(z)$ at the origin,
\begin{equation}\label{RIntDef}
      {}^R\int_M \mu = \FP_{z=0} \zeta_x(z).
\end{equation}

Alternately, we could consider the integral of $\mu$ on the truncated manifold $\{x \geq \varepsilon \}$. The phg expansion of $\mu$ induces a phg expansion of
\begin{equation*}
      \int_{x \geq \varepsilon} \mu
\end{equation*}
in $\varepsilon$, and we could define a renormalized integral of $\mu$ as the coefficient of $\varepsilon^0$. 

Both of these renormalization schemes are widely used in the literature. The latter is the definition used in \cite{APS Book} to define the ``$b$-integral", for instance. We refer the reader to \cite{My Preprint} for a  comparison of the two schemes. For the integrals of interest here, these two definitions coincide, and we shall mostly use \eqref{RIntDef}.

%%%%%%%%%%%%%%%%%%%%%%%%%%%%%%%%%%%%%%%%%%%%%%%%%%%%%%%%%%%%%%%%%%%%%%%%%%%%%%%%%%%%%%%%%%%%%%%%%%%%%%
\subsubsection{Renormalized Dimension} $ $\newline
%%%%%%%%%%%%%%%%%%%%%%%%%%%%%%%%%%%%%%%%%%%%%%%%%%%%%%%%%%%%%%%%%%%%%%%%%%%%%%%%%%%%%%%%%%%%%%%%%%%%%%

If $A$ is a trace-class operator that acts through an integral kernel
\begin{equation*}
      Af(\xi) = \int \curly{K}_A\lrpar{\xi, \xi'} f\lrpar{\xi'} \; \mathrm{d}\xi',
\end{equation*}
then Lidskii's theorem expresses its trace as the integral of its kernel along the diagonal,
\begin{equation*}
      \mathrm{Tr}(A) = \int \curly{K}_A\lrpar{ \xi, \xi} \; \mathrm{d}\xi.
\end{equation*}
This applies in particular to smoothing pseudo-differential operators on closed manifolds.

A smoothing operator on an edge manifold is also given by integrating against a smooth kernel, but its restriction to the diagonal fails to be integrable. Nevertheless, kernels of edge pseudo-differential operators have, by definition, phg expansions at the boundary faces of the double edge space. In particular, the restriction to the diagonal of a smoothing operator will have a phg expansion at the boundary and we can use the discussion from \ref{subsec:RInt} to define its renormalized integral, hereafter known as its renormalized trace,
\begin{equation*}
      {}^R\mathrm{Tr}(A) = {}^R\int \curly{K}_A\lrpar{ \xi, \xi} \; \mathrm{d}\xi .
\end{equation*}

Finally, assume that the projection, $\curly{P}$ onto a subspace of $L^2(M)$ is an element of $\Psi^{-\infty, \curly{E}}_e$. We define the renormalized dimension of this space by
\begin{equation*}
      {}^R\mathrm{dim} = {}^R\mathrm{Tr}\lrpar{\curly{P}}.
\end{equation*}
This applies to the spaces of harmonic forms, or more generally, the null spaces of elliptic edge pseudo-differential operators by the results of \cite{Mazzeo:Edge}.

Note that the renormalized dimension is only a dimension by analogy. It is {\em a priori} neither positive nor an integer. Nevertheless, we shall see that it comes up naturally in the heat equation proof of the index theorem.

%%%%%%%%%%%%%%%%%%%%%%%%%%%%%%%%%%%%%%%%%%%%%%%%%%%%%%%%%%%%%%%%%%%%%%%%%%%%%%%%%%%%%%%%%%%%%%%%%%%%%%
\subsection{The Index Theorem} $ $\newline
%%%%%%%%%%%%%%%%%%%%%%%%%%%%%%%%%%%%%%%%%%%%%%%%%%%%%%%%%%%%%%%%%%%%%%%%%%%%%%%%%%%%%%%%%%%%%%%%%%%%%%

In analogy to McKean-Singer, consider the identity
\begin{equation}\label{McSing}
      \lim_{t \to \infty} {}^R\mathrm{Str}\lrpar{ e^{-t\eth^2}} -
      \lim_{t \to 0} {}^R\mathrm{Str}\lrpar{ e^{-t\eth^2}} =
      \int_{0}^{\infty} \partial_t {}^R\mathrm{Str}\lrpar{ e^{-t\eth^2}} \; \mathrm{dt}.
\end{equation}
Note that both sides depend on the choice of bdf, $x$.

\begin{lemma}
If the image of $\eth^2$ is closed, then
\begin{equation}\label{t=infty}
      \lim_{t \to \infty} {}^R\mathrm{Str}\lrpar{ e^{-t\eth^2}} =
       {}^R\mathrm{Ind}\lrpar{ \eth } .
\end{equation}
\end{lemma}

{\em Remark.} By the final remark of section $\mathcal{x}$\ref{sec:elliptic}, this lemma applies to the Laplacian of an exact edge metric with $\dim B >0$.

\begin{proof}
To prove the lemma, it suffices to show that the heat kernel converges uniformly to the integral kernel of the projection along the diagonal.
Let $\curly{P}$ be the projection onto $\ker \eth$, and consider $\eth_0 = \eth - \curly{P}$. As $\eth_0^2$ is a positive injective operator with closed image its spectrum has a positive lower bound, $\lambda_1$. 
Thus, for any $\omega$,
\begin{equation*}
      0 \leq \left< e^{-t\eth_0^2}\omega, \omega \right> \leq e^{-t\lambda_1} \left< \omega, \omega \right>.
\end{equation*}
This $L^2$-control on the heat kernel translates, via its semigroup property, to uniform pointwise control.
Indeed, if we denote the integral kernel of $e^{-t\eth^2_0}$ by $\curly{K}_t$, then from
\begin{equation*}
	e^{-\lrpar{s+t}\eth_0^2}=e^{-s\eth_0^2}e^{-t\eth_0^2} \iff
	\curly{K}_{s+t}\lrpar{\zeta, \zeta'} = \int_{\zeta''} \curly{K}_s\lrpar{\zeta, \zeta''}
		\curly{K}_t\lrpar{\zeta'', \zeta'},
\end{equation*}
we see that
\begin{equation*}\begin{split}
	\curly{K}_{t+2}\lrpar{\zeta,\zeta} &=
	\int_{\zeta', \zeta''} \curly{K}_{1}\lrpar{\zeta, \zeta'} \curly{K}_t\lrpar{\zeta',\zeta''}
		\curly{K}_1\lrpar{\zeta'', \zeta}
	= \left< e^{t\eth_0^2}\curly{K}_1\lrpar{\zeta, \cdot}, \curly{K}_1\lrpar{\zeta, \cdot} \right>\\
	&\leq e^{-t\lambda_1}\left<\curly{K}_1\lrpar{\zeta, \cdot}, \curly{K}_1\lrpar{\zeta, \cdot}\right>
	= e^{-t\lambda_1}\curly{K}_2\lrpar{\zeta, \zeta}.
\end{split}\end{equation*}
By Theorem \ref{HeatKerExt}, $\curly{K}_2$ restricts to the diagonal as a smooth function on a compact set, so we have uniform convergence of $K_t\lrpar{\zeta,\zeta}$ to zero.
The lemma follows from the relation
\begin{equation*}
	e^{-t\eth_0^2} = \lrpar{e^{-t}-1}\curly{P}+e^{-t\eth^2}.
\end{equation*}

\end{proof}

The work for the limit as $t \to 0$ has already been done in \cite{APS Book}. Indeed, by the local index theorem we know that the (pointwise) supertrace of the heat kernel in the interior of $M$ tends to the Atiyah-Singer integrand, $\mathrm{AS}$,
\begin{equation*}
      \mathrm{str}\lrpar{e^{-t\eth^2}}\lrpar{\zeta} \to \mathrm{AS}\lrpar{\zeta},
\end{equation*}
for any $\zeta \in M$. Indeed, the convergence is uniform in $C^{\infty}$.
As discussed in \cite{APS Book},
\begin{equation}\label{t=0}
      \lim_{t \to 0} {}^R\mathrm{Str}\lrpar{ e^{-t\eth^2}} = {}^R\int_M \mathrm{AS} 
\end{equation}
can be thought of as following by continuity.

As for the final term in \eqref{McSing}, note that
\begin{equation}\label{GenMcSing}\begin{split}
      \partial_t \mathrm{Str} \lrpar{ x^z e^{-t\eth^2} } 
      &= -\frac{1}{2} \mathrm{Str} \lrpar{x^z \lrspar{\eth, \eth e^{-t\eth^2}} } 
      =  -\frac{1}{2} \mathrm{Str} \lrpar{ \lrspar{x^z, \eth}\eth e^{-t\eth^2} %
      		+ \lrspar{\eth, x^z\eth e^{-t\eth^2}} }\\
      &=  -\frac{1}{2} \mathrm{Str} \lrpar{ \lrspar{x^z, \eth}\eth e^{-t\eth^2} }
      =  \frac{1}{2} \mathrm{Str} \lrpar{ \cl\lrpar{d\lrpar{x^z}}\eth e^{-t\eth^2} }
      = \frac{z}{2} \mathrm{Str} \lrpar{ x^z \cl\lrpar{ \frac{dx}{x} } \eth e^{-t\eth^2} },
\end{split} \end{equation}
hence,
\begin{equation}\label{Reta}
      2\int_{0}^{\infty} \partial_t {}^R\mathrm{Str}\lrpar{ e^{-t\eth^2}} \; \mathrm{dt}
      =  \int_{0}^{\infty} 
            \int_{\partial M} \lrspar{ \mathrm{str}%
               \lrpar{ \cl\lrpar{ \frac{dx}{x} } \eth e^{-t\eth^2}} %
               \downharpoonright_{\mathrm{diag}} }_{(-1)}
            \; \mathrm{dt}.
\end{equation}
This is similar to one of the standard definitions of the $\eta$ invariant. Indeed, in the asymptotically cylindrical case, this is the usual $\eta$ invariant of the boundary, see \cite{APS Book}. We will refer to this as the ``renormalized $\eta$ invariant" and denote it by ${}^R\eta$. In the next section, we will show that for the Gauss-Bonnet complex, ${}^R\eta =0$. A more detailed analysis of ${}^R\eta$ will be carried out in \cite{Me-Rafe}, via the ``Getzler rescaling" technique of \cite[Ch. 8]{APS Book}.

We put \eqref{t=infty}, \eqref{t=0}, and \eqref{Reta} into \eqref{McSing} and obtain the following theorem.
\begin{theorem}\label{IndThm}
If $\eth$ is a generalized edge Dirac operator, then 
\begin{equation*}
      {}^R\int \mathrm{AS} - \frac{1}{2}{}^R\eta = \lim_{t \to \infty} {}^R\mathrm{Str}\lrpar{\eth}.
\end{equation*}
If the image of $\eth$ is closed, then this limit is ${}^R\mathrm{Ind}\lrpar{\eth}$.
\end{theorem}
Finally, except for the construction of the heat kernel, we have used nothing specific to edge metrics not shared by other MICE and this theorem is true for any such metric once the heat kernel has been constructed. On the other hand, Theorem \ref{IndThm} is unsatisfactory in that two of the three terms remain mysterious. In the next section we will remedy this for the Gauss-Bonnet operator.

%%%%%%%%%%%%%%%%%%%%%%%%%%%%%%%%%%%%%%%%%%%%%%%%%%%%%%%%%%%%%%%%%%%%%%%%%%%%%%%%%%%%%%%%%%%%%%%%%%%%%%
\section{The Gauss-Bonnet Theorem} \label{sec:GBThm}
%%%%%%%%%%%%%%%%%%%%%%%%%%%%%%%%%%%%%%%%%%%%%%%%%%%%%%%%%%%%%%%%%%%%%%%%%%%%%%%%%%%%%%%%%%%%%%%%%%%%%%

Consider the de Rham operator, 
\begin{equation*}
      \Omega^*\lrpar{M} \xrightarrow{d+\delta} \Omega^*\lrpar{M}.
\end{equation*}
The space of forms on an even dimensional manifold is naturally a $\mathbb{Z}/2$-graded Clifford module with respect to the splitting
\begin{equation*}
      \Omega^*\lrpar{M} = \Omega^{\mathrm{even}} \lrpar{M} \oplus \Omega^{\mathrm{odd}} \lrpar{M},
\end{equation*}
and the Clifford action
\begin{equation*}
      \cl\lrpar{ \theta} \omega = \theta \wedge \omega - \theta \hook \omega .
\end{equation*}
The Gauss-Bonnet operator, $\eth_{\mathrm{GB}}$, is the corresponding Dirac operator.
As we have anticipated, in this situation the renormalized $\eta$ invariant is trivial.

\begin{theorem}
For the Gauss-Bonnet complex,
\begin{equation*}
	\mathrm{str} \lrpar{ \cl\lrpar{ \frac{dx}{x} } \eth e^{-t\eth^2}} =0 ,
\end{equation*}
hence, from \eqref{GenMcSing}, $\mathrm{Str} \lrpar{ x^z e^{-t\eth^2} } $ is independent of $t$.
In particular, this implies that ${}^R\eta\lrpar{\eth_{GB}} = 0$.
\end{theorem}

\begin{proof}
Note that
\begin{equation*}
      	\Omega_{x}^{even}(\partial M) \oplus \frac{dx}{x}\wedge\Omega_{x}^{odd}(\partial M)
	     \xrightarrow{\eth_{GB}} 
	\Omega_{x}^{odd}(\partial M) \oplus \frac{dx}{x}\wedge\Omega_{x}^{even}(\partial M)
\end{equation*}
is given by
\begin{equation*}
        \begin{pmatrix}
	    \widehat{\eth_{GB}} & -\nabla_{x\partial_{x}} \\
	    \nabla_{x\partial_{x}} & -\widehat{\eth_{GB}}
	\end{pmatrix}
\end{equation*}
where $\widehat{\eth_{GB}} = \eth_{GB} - \cl(\frac{dx}{x})\nabla_{x\partial_{x}}$.

Similarly, with respect to this splitting, we have:
\begin{equation*}
	\Delta = \begin{pmatrix}
	    \widehat{\eth_{GB}}^{2}-\nabla^{2} & [\nabla,\widehat{\eth_{GB}}] \\
	    [\nabla,\widehat{\eth_{GB}}] & \widehat{\eth_{GB}}^{2}-\nabla^{2}
	\end{pmatrix}, \phantom{xxx}
	\cl(\frac{dx}{x}) = \begin{pmatrix}
	    0 & -1 \\
	    1 & 0
	\end{pmatrix},
\end{equation*}
hence
\begin{equation*}
       \cl(\frac{dx}{x})\eth_{GB}e^{-t\Delta} = \begin{pmatrix}
           -\nabla_{x\partial_{x}} & \widehat{\eth_{GB}} \\
           \widehat{\eth_{GB}} & -\nabla_{x\partial_{x}}
        \end{pmatrix} \exp \left( -t \begin{pmatrix}
           \widehat{\eth_{GB}}^{2}-\nabla^{2} & [\nabla,\widehat{\eth_{GB}}] \\
           [\nabla,\widehat{\eth_{GB}}] & \widehat{\eth_{GB}}^{2}-\nabla^{2}
        \end{pmatrix} \right)
        =: \begin{pmatrix}
           A_{ee} & B_{oe} \\
           B_{eo} & A_{oo}
        \end{pmatrix}.
\end{equation*}
Note that this is formally the same expression if we were to interchange the even and odd parts,
\begin{equation*}
     \Omega_{x}^{odd}(\partial M) \oplus \frac{dx}{x}\wedge\Omega_{x}^{even}(\partial M)
		\rightarrow 
     \Omega_{x}^{odd}(\partial M) \oplus \frac{dx}{x}\wedge\Omega_{x}^{even}(\partial M).
\end{equation*}
That is, we have
\begin{equation*}
      \Omega^{e}\oplus\Omega^{o}
	\xrightarrow{ \begin{pmatrix}
    		A_{ee} & B_{oe} \\ B_{eo} & A_{oo} \end{pmatrix} }
      \Omega^{e}\oplus\Omega^{o} \phantom{xxx}
      \Omega^{o}\oplus\Omega^{e}
	\xrightarrow{ \begin{pmatrix}
		A_{oo} & B_{eo} \\ B_{oe} & A_{ee} \end{pmatrix} }
      \Omega^{o}\oplus\Omega^{e},
\end{equation*}
so in either case the trace is
\begin{equation*}
      \mathrm{tr} \left( \Omega^{e} \xrightarrow{A_{ee}} \Omega^{e} \right) 
	+ \mathrm{tr} \left( \Omega^{o} \xrightarrow{A_{oo}} \Omega^{o} \right),
\end{equation*}
and the supertrace vanishes.
\end{proof}

Thus we know from Theorem \ref{IndThm} that
\begin{equation}\label{RGB}
      {}^R\int \mathrm{Pff} = \lim_{t\to \infty} {}^R\mathrm{Str}\lrpar{e^{-t\Delta}}.
\end{equation}
We can supplement this equation through Chern's Gauss-Bonnet theorem for incomplete metrics on manifolds with boundary. We think of the resulting limiting formula as a ``soft" index theorem in contrast to \eqref{RGB}.

%%%%%%%%%%%%%%%%%%%%%%%%%%%%%%%%%%%%%%%%%%%%%%%%%%%%%%%%%%%%%%%%%%%%%%%%%%%%%%%%%%%%%%%%%%%%%%%%%%%%%%
\subsection{Soft Index Formula} $ $\newline
%%%%%%%%%%%%%%%%%%%%%%%%%%%%%%%%%%%%%%%%%%%%%%%%%%%%%%%%%%%%%%%%%%%%%%%%%%%%%%%%%%%%%%%%%%%%%%%%%%%%%%

On an even-dimensional manifold with boundary, $(X,g)$, Chern's Gauss-Bonnet index formula is
\begin{equation}\label{ChernGB}
      \int_X \mathrm{Pff} + \int_{\partial X} \Fun = \chi(X),
\end{equation}
where $\mathrm{Pff}$ is the Pfaffian and $\Fun$ is a polynomial in the curvature and the second fundamental form. This easily yields the following formula.

\begin{theorem}\label{ChernRGB}
Let $(M,g)$ be the interior of a compact manifold with boundary, $x$ a bdf, and assume that $g$ has a phg expansion at $\partial M$ in terms of $x$. Then the Pfaffian and $\ITFun$ also have phg expansions in $x$, so that we may define their renormalized integrals as in section $\mathcal{x}$\ref{subsec:RInt}. Furthermore, these satisfy
\begin{equation}
      {}^R\int \mathrm{Pff} + \FP_{\varepsilon =0} \int_{x=\varepsilon} \ITFun = \chi (M).
\end{equation}
\end{theorem}

\begin{proof}
The Pfaffian and $\Fun$ are polynomials in the curvature and the second fundamental form, hence inherit a phg expansion from that of $g$.
Now, we simply consider \eqref{ChernGB} for the manifold $\{ x \geq \varepsilon \}$,
\begin{equation*}
      \int_{\{ x \geq \varepsilon \}} \mathrm{Pff} + \int_{x = \varepsilon} \Fun 
         = \chi \lrpar{\{ x \geq \varepsilon \}}.
\end{equation*}
For small enough $\varepsilon$ the right hand side is $\chi(M)$, hence independent of $\varepsilon$. Thus the left hand side must be independent of $\varepsilon$. So the equality is true replacing the left hand side with the $\varepsilon^{0}$ term in its expansion, and this proves the theorem.
\end{proof}

For the next section it will be useful to have explicit expressions for the Pfaffian and $\Fun$ on a manifold of even dimension $m=2n$. For an arbitrary local frame, denote the curvature by $R$ and the second fundamental form by $\SFun$, then the Pfaffian is given by
\[       \mathrm{Pff} = \frac{(-1)^n}{2^{3n}\pi^nn!} \sum_{\sigma, \tau \in \Sigma_{2n}} (-1)^{|\sigma| + |\tau|}
         R^{\sigma_1\sigma_2}_{\tau_1\tau_2} \cdots 
         R^{\sigma_{2n-1}\sigma_{2n}}_{\tau_{2n-1}\tau_{2n}}  \dvol,\]
and similarly
\begin{equation*}
     \Fun = \sum_{q=0}^{n-1}  %
         \frac{(-1)^q}{2^{3q}\pi^qq! \mathrm{vol}\lrpar{\mathbb{S}^{m-1-2q}}\lrpar{m-1-2q}!}
         Q_{q.m} \dvol_{\partial M}
\end{equation*}
with
\begin{equation*}
         Q_{q,m} = \sum_{\sigma, \tau \in \Sigma_{m-1}} (-1)^{|\sigma| + |\tau|}%
         R^{\sigma_1\sigma_2}_{\tau_1\tau_2} \cdots %
         R^{\sigma_{2q-1}\sigma_{2q}}_{\tau_{2q-1}\tau_{2q}} % 
         \SFun^{\sigma_{2q+1}}_{\tau_{2q+1}} \cdots %
         \SFun^{\sigma_{m-1}}_{\tau_{m-1}} . 
\end{equation*}
We will use the formalism of double forms as set out in \cite{Kulkarni} and \cite{Tubes} to compute with the curvature and the second fundamental form, see \cite[$\mathcal{x}$4]{My Preprint}.

%%%%%%%%%%%%%%%%%%%%%%%%%%%%%%%%%%%%%%%%%%%%%%%%%%%%%%%%%%%%%%%%%%%%%%%%%%%%%%%%%%%%%%%%%%%%%%%%%%%%%%
\subsection{Fibrations with Trivial Base or Fiber} \label{sec:TrivFib} $ $\newline
%%%%%%%%%%%%%%%%%%%%%%%%%%%%%%%%%%%%%%%%%%%%%%%%%%%%%%%%%%%%%%%%%%%%%%%%%%%%%%%%%%%%%%%%%%%%%%%%%%%%%%

The simplest boundary fibration structures correspond to fibrations with zero-dimensional bases or fibers.  These include asymptotically cylindrical metrics (b-calculus), conformally compact and asymptotically hyperbolic metrics ($0$-calculus), as well as asymptotically flat (scattering calculus) and others. Near the boundary, these metrics can be put in the form
\begin{equation}\label{GenMetric}
	 \frac{dx^2}{\alpha^2x^{2\eta}} + \frac{h_x}{x^{2\beta}},
\end{equation}
where $\eta \geq 1$ and $\beta \geq 0$ are constants, $\eta \geq \beta$, $h_x$ is a family of metrics on $\partial M$, and $\alpha$ is the pointwise length of $|dx|$ with respect to the metric $\bar{g}$ below -- in particular it does not vanish when $x=0$. 
The analysis of $\Fun$ for these metrics is similar in the appearance of the Weyl volume of tubes invariants for $(\partial M, h_x)$. Nevertheless, only in the simplest situations ($\beta=0$ or $\beta=\eta -1$) can we figure out $\FP_{\varepsilon =0} \int_{x=\varepsilon} \Fun$ by direct computation.

Our study of $\Fun$ proceeds by comparing the curvature, $R$, and second fundamental form, $\SFun$, of $(M,g)$ with the corresponding tensors for 
\begin{equation*}
	\bar{g}:= \frac{dx^2}{\alpha^2} + h_x.
\end{equation*}
Tensors corresponding to $\bar{g}$ will be differentiated from those corresponding to $g$ by the presence of a bar.
Since $\Fun$ only involves directions tangent to the boundary, we can use the Gauss equation to compare the curvature $R$ of $g$ and $\bar{R}$ of $\bar{g}$ by comparing $R^{\partial M}$ and $\bar{R}^{\partial M}$ as well as $\SFun$ for $g$ and $\bar{g}$. Note that
\begin{equation*}
	R^{\partial M}=x^{2\beta}\bar{R}^{\partial M}.
\end{equation*}

For the second fundamental form, consider a local frame (Fermi coordinates), $\{ \bar{X}_i \}$, for $\bar{g}$ centered at a point $p \in \partial M$ with $\bar{X}_i$ orthogonal to $\partial_x=:\bar{X}_m$, and the corresponding local frame, $X_i:=x^{\beta}\bar{X}_i$, $X_m = x^{\eta}\bar{X}_m$ for $g$ (cf. \cite[$\mathcal{x}$3]{My Preprint}).
The second fundamental forms are given by
\begin{equation*}
	\bar{\SFun}\lrpar{\bar{X}_i,\bar{X}_j}
	= \bar{g} \lrpar{ \bar{\nabla}_{\bar{X}_i} \bar{X}_j, \alpha\bar{X}_m }
	= \frac{1}{\alpha} \bar{\Gamma}_{ij}^m, \phantom{xxx}
	\SFun\lrpar{X_i, X_j} = \frac{1}{\alpha} \gamma_{ij}^m.
\end{equation*}
These are easily compared using the Koszul formula
\begin{equation*}\begin{split}
	\SFun\lrpar{X_i,X_j} = \alpha g\lrpar{ \nabla_{X_i} X_j, X_m } %
	&= \frac{\alpha}{2} 
	   \left[ X_i g\lrpar{ X_j, X_m } + X_j g\lrpar{ X_m, X_i } - X_m g\lrpar{ X_i, X_j} \right.\\
	&\left. \phantom{xxxxx} + g\lrpar{ \lrspar{X_i , X_j}, X_m } + g\lrpar{ \lrspar{X_m, X_i}, X_j} %
	   - g\lrpar{ \lrspar{X_j, X_m}, X_i} \right] \\
	&= \frac{\alpha}{2} \left[ - X_m g\lrpar{ X_i, X_j}  +  g\lrpar{ \lrspar{X_m, X_i}, X_j} %
	   + g\lrpar{ \lrspar{X_m, X_j}, X_i} \right] \\
	&= x^{\eta}\bar{\SFun}\lrpar{\bar{X}_i,\bar{X}_j}  + \alpha \beta x^{\eta-1} \bar{g}_{ij}.
\end{split}\end{equation*}
Hence as double forms,
\begin{equation}\label{SecBar}
	\FunForm\lrpar{X_i,X_j} 
	= x^{\eta-1}\lrpar{x\bar{\FunForm} + \alpha \beta \bar{\df{g}}}
	   \lrpar{\bar{X}_i, \bar{X}_j}.
\end{equation}

Using the Gauss equation, we have found the change in curvature (for vector fields tangent to the level sets of $x$)
\begin{equation}\label{GenCase}\begin{split}
	\df{R}\lrpar{X_i,X_j}\lrpar{X_k, X_{\ell}} 
	&= \lrpar{ \df{R}^{\partial M} - \frac{\FunForm^2}{2}}\lrpar{X_i,X_j}\lrpar{X_k, X_{\ell}} \\
	&= \lrpar{ x^{2\beta}\bar{\df{R}}^{\partial M} %
	   - \frac{\lrpar{x^{\eta-1}\lrpar{x\bar{\FunForm} + %
	      \alpha \beta \bar{\df{g}}}}^2}{2}}\lrpar{\bar{X}_i,\bar{X}_j}\lrpar{\bar{X}_k,\bar{X}_{\ell}} \\
	&= \lrpar{ x^{2\beta}\bar{\df{R}}^{\partial M} %
	   - x^{2\eta} \frac{\bar{\FunForm}^2}{2} %
	   - x^{2\eta-1}\alpha\beta \bar{\FunForm}\bar{g} %
	   - x^{2\eta-2}\alpha^2\beta^2 \frac{\bar{g}^2}{2} }
	   \lrpar{\bar{X}_i,\bar{X}_j}\lrpar{\bar{X}_k, \bar{X}_{\ell}}.
\end{split}\end{equation}
Note the different behaviors for different values of $\eta$, $\beta$. For asymptotically cylindrical ends ($\eta=1$, $\beta=0$),
\begin{equation*}
	{}^b\df{R}\lrpar{X_i,X_j}\lrpar{X_k, X_{\ell}} 
	= \lrpar{ \bar{\df{R}}^{\partial M} %
	   - x^{2} \frac{\bar{\FunForm}^2}{2} }
	   \lrpar{\bar{X}_i,\bar{X}_j}\lrpar{\bar{X}_k, \bar{X}_{\ell}}.
\end{equation*}
On the other hand conformally compact metrics ($\eta = \beta =1$) are asymptotically isotropic,
\begin{equation*}
	{}^0\df{R}\lrpar{X_i,X_j}\lrpar{X_k, X_{\ell}} 
	= \lrpar{ x^{2}\bar{\df{R}}^{\partial M} %
	   - x^{2} \frac{\bar{\FunForm}^2}{2} %
	   - x\alpha \bar{\FunForm}\bar{g} %
	   - \alpha^2\frac{\bar{g}^2}{2} }
	   \lrpar{\bar{X}_i,\bar{X}_j}\lrpar{\bar{X}_k, \bar{X}_{\ell}}.
\end{equation*}
In particular, since the left hand side does not depend on the choice of $x$, neither does $\alpha \downharpoonright_{\partial M}$. In fact, from \eqref{GenCase}, this is true whenever $\beta \notin \{0,\eta-1\}$. A bdf for a conformally compact metric is called special if $|dx| = \alpha \downharpoonright_{\partial M}$ on a neighborhood of the boundary. Given any bdf $x_0$ and a conformally compact metric $g$, there exists a special bdf $x$ with 
$x_0^2g \downharpoonright_{\partial M} = x^2g \downharpoonright_{\partial M}$
Conformally compact manifolds with $\alpha\downharpoonright_{\partial M} \equiv 1$ are known as asymptotically hyperbolic. These include the Poincar\'e-Einstein manifolds of the AdS/CFT correspondence in physics.
Another interesting particular case is that of scattering metrics ($\eta=2$, $\beta =1$) which are asymptotically flat,
\begin{equation*}
	{}^{\mathrm{sc}}\df{R}\lrpar{X_i,X_j}\lrpar{X_k, X_{\ell}} 
	= \lrpar{ x^{2}\bar{\df{R}}^{\partial M} %
	   - x^{4} \frac{\bar{\FunForm}^2}{2} %
	   - x^{3}\alpha \bar{\FunForm}\bar{g} %
	   - x^{2}\alpha^2 \frac{\bar{g}^2}{2} }
	   \lrpar{\bar{X}_i,\bar{X}_j}\lrpar{\bar{X}_k, \bar{X}_{\ell}}.	
\end{equation*}
Radial compactification of $\mathbb{R}^n$ to a half-sphere of the same dimension produces a scattering metric.

What does \eqref{GenCase} tell us about $\Fun$? We can alternately think of $\df{R}$ as
\begin{equation*}
	\Omega^2 \xrightarrow{\df{R}} \Omega^2
	\text{ or } \df{R} \in \Omega^2 \otimes \Omega^2,
\end{equation*}
with the latter yielding the coefficients for the former when viewed as a two-form. We have found a relation between $\df{R}(X_i,X_j)(X_k,X_{\ell})$ and some other double form, say $\df{S}$, evaluated at $(\bar{X}_i, \bar{X}_j)(\bar{X}_k, \bar{X}_{\ell})$, i.e. we have expressions for the coefficients of $\df{R}$. This means that
\begin{equation*}
	\left< \{ X_s \} \right> \xrightarrow{\df{R}} \left< \{ X_s \} \right> \iff
	\left< \{ \bar{X}_s \} \right> \xrightarrow{\frac{1}{x^{\beta(n-1)}}\df{S}} \left< \{ \bar{X}_s \} \right>.
\end{equation*}
Also note that both \eqref{SecBar} and \eqref{GenCase} are polynomials in
\begin{equation*}
	x^{2\beta}\bar{R}^{\partial M}, x^{\eta}\bar{\FunForm}, \text { and } x^{\eta-1}\alpha\beta\bar{g}.
\end{equation*}
Thus we can conclude that, for some constants $C_{a,b,c}$,
\begin{equation}\label{FunCon}
	\Fun = x^{-\beta\lrpar{n-1}}
	\sum_{2a+b+c=n-1} C_{a,b,c} x^{2a\beta+b\eta + c(\eta-1)} \alpha^c %
	\lrspar{ \lrpar{\bar{\df{R}}^{\partial M}}^a \bar{\FunForm}^b \bar{g}^c %
	\lrpar{\bar{X}_1,\ldots \bar{X}_{n-1}} }.
\end{equation}

This expression simplifies when $\beta =0$,
\begin{equation}\label{beta=0}
	\Fun =
	\sum_{2a+b=n-1} C_{a,b} x^{b\eta} 
	\lrspar{ \lrpar{\bar{\df{R}}^{\partial M}}^a \bar{\FunForm}^b %
	\lrpar{\bar{X}_1,\ldots \bar{X}_{m-1}} }	
\end{equation}
and when $\beta = \eta-1$,
\begin{equation}\label{beta=eta-1}
	\Fun =
	\sum_{2a+b+c=n-1} C_{a,b,c} x^{b} \alpha^c
	\lrspar{ \lrpar{\bar{\df{R}}^{\partial M}}^a \bar{\FunForm}^b \bar{g}^c %
	\lrpar{\bar{X}_1,\ldots \bar{X}_{m-1}} }.
\end{equation}
Thus we have shown the following consequence of the soft index theorem (Theorem \ref{ChernRGB}).

\begin{corollary}\label{B/SC}
Assume that $M$ is an even-dimensional manifold with a metric of the form \eqref{GenMetric}.
If $\beta=0$, e.g. for a asymptotically cylindrical metric or a cusp metric,
\begin{equation}\label{bRGB}
	\FP_{\varepsilon =0} \int_{x=\varepsilon} \ITFun = 0, \text{ hence }
	{}^R\int \mathrm{Pff} = \chi(M).
\end{equation}
If instead we have $\beta = \eta -1$ and $\alpha \equiv 1$, such as for a scattering metric, then
\begin{equation}\label{scRGB}
	\FP_{\varepsilon =0} \int_{x=\varepsilon} \ITFun = P_{\eta,n}\lrpar{\partial M, h_0}
\end{equation}
is a linear combination of the Weyl volume of tubes invariants of the boundary metric $h_0$, hence
\begin{equation*}
	{}^R\int \mathrm{Pff} + P_{\eta,n}\lrpar{\partial M, h_0} = \chi(M).	
\end{equation*}
\end{corollary}

\begin{proof}
The proof of \eqref{bRGB} follows directly from \eqref{beta=0} which shows that $\Fun$ vanishes with $x$.
Similarly, when $\beta = \eta -1$, \eqref{beta=eta-1} shows that
\begin{equation}\begin{split}
	\FP_{\varepsilon =0} \int_{x=\varepsilon} \Fun
	&= \sum_{2a+c=n-1} C_{a,0,c} \alpha^c
		\lrspar{ \lrpar{\bar{\df{R}}^{\partial M}}^a \bar{g}^c %
		\lrpar{\bar{X}_1,\ldots \bar{X}_{m-1}} }\\
	&= \frac{1}{\lrpar{2\pi}^{n/2}}\sum_{q=0}^{n/2-1}\sum_{j=0}^{q}
		\frac{(-1)^j\beta^{n-1-2q}\alpha^{n-1-2j}}{\lrpar{n-1-2q}!!j!\lrpar{q-j}!2^{q-j}}
		\lrpar{\bar{\df{R}}^{\partial M}}^j \bar{\df{g}}^{n-1-2j}
		\lrpar{\bar{X}_1,\ldots \bar{X}_{m-1}} 
.\end{split}\end{equation}
\end{proof}

%%%%%%%%%%%%%%%%%%%%%%%%%%%%%%%%%%%%%%%%%%%%%%%%%%%%%%%%%%%%%%%%%%%%%%%%%%%%%%%%%%%%%%%%%%%%%%%%%%%%%%
\section{Finite Time Trace on Conformally Compact Manifolds} \label{sec:ConfComp}
%%%%%%%%%%%%%%%%%%%%%%%%%%%%%%%%%%%%%%%%%%%%%%%%%%%%%%%%%%%%%%%%%%%%%%%%%%%%%%%%%%%%%%%%%%%%%%%%%%%%%%

A particular type of edge metric, when the fibers of the boundary fibration consist of a point, is known as a conformally compact metric. So called because the metric $\bar{g}=x^2g$, where $x$ is any bdf, extends to a metric on $\overline{M}$.
In terms of \eqref{GenMetric}, conformally compact metrics correspond to $\beta=\eta=1$. As we remarked above, for these metrics it is possible to choose a ``special" bdf, $x$, so that $\alpha \equiv \alpha\downharpoonright_{\partial M}$. Throughout this section, $x$ will always denote a special bdf.

The study of these metrics began in \cite{Mazzeo:Hodge} and \cite{Mazzeo-Melrose}, where their Hodge theory and the meromorphic continuation of the resolvent were carried out, respectively. In particular, from \cite{Mazzeo:Hodge} we know that $\eth_{GB}$ has closed image and that the spaces of $L^2$ harmonic forms have topological interpretation:
\begin{equation*}
      \mathcal{H}_{L^2}^k = \begin{cases}
         H^k(M, \partial M) & \text{ if $k<\frac{m}{2}$} \\
         H^k(M) & \text{ if $k>\frac{m}{2}$}
      \end{cases},
\end{equation*}
where we are assuming that $M$ is even dimensional. The $L^2$ harmonic forms in middle degree form an infinite dimensional space, essentially because of the conformal invariance of the $L^2$-norm in middle degree. Thus the Gauss-Bonnet Theorem \ref{RGB} in this context is

\begin{corollary} \label{ConfComp1}
For any choice of bdf, $x$, on an even-dimensional conformally compact manifold,
\begin{equation*}\begin{split}
      {}^R\int \mathrm{Pff} &=
         \sum_{k<\frac{m}{2}} (-1)^k \dim H^k(M,\partial M)
         +(-1)^{m/2}  \lrpar{ {}^R\dim \mathcal{H}_{L^2}^{m/2}}
         +\sum_{k>\frac{m}{2}} (-1)^k \dim H^k(M)\\
         &= 2\sum_{k<\frac{m}{2}} (-1)^k \dim H^k(M,\partial M)
	+ (-1)^{m/2}  \lrpar{ {}^R\dim \mathcal{H}_{L^2}^{m/2}}.
\end{split} \end{equation*}
\end{corollary}
Note that both sides depend on the choice of $x$, though the dependence of the renormalized integral of the Pfaffian is exactly compensated by that of the renormalized dimension of middle-degree harmonic forms.

In this context, \eqref{FunCon} becomes
\begin{equation*}
	\Fun =
	\sum_{2a+b+c=n-1} C_{a,b,c} \frac{\alpha^c}{x^c} %
	\lrspar{ \lrpar{\bar{\df{R}}^{\partial M}}^a \bar{\FunForm}^b \bar{g}^c %
	\lrpar{\bar{X}_1,\ldots \bar{X}_{n-1}} }.
\end{equation*}
Thus the constant term in the expansion of $\int_{x=\varepsilon} \Fun$ involves integrating over the boundary the $x^c$ term in the expansion of $\lrspar{ \lrpar{\bar{\df{R}}^{\partial M}}^a \bar{\FunForm}^b \bar{g}^c \lrpar{\bar{X}_1,\ldots \bar{X}_{n-1}} }$, and in general there does not seem to be a simple approach to understanding these terms.

There are special classes of conformally compact manifolds for which we can show that 
\begin{equation*}
	\FP_{\varepsilon =0} \int_{x=\varepsilon} \Fun =0.
\end{equation*}
This is true whenever the family of tensors $h_x$ from \eqref{GenMetric} has an expansion in $x$ involving only {\em even} powers of $x$ below $x^n$. It is a nice property of these manifolds that this property is independent of the choice of special bdf \cite{Guillarmou}. This can be traced back to the following very useful fact about special bdfs. If $\hat{x}=e^{\omega(x,y)} x$ and $x$ are both special bdfs, and $h_x$ has only even powers of $x$ in its expansion below $x^\ell$, then $\omega(x,y)$ has only even powers of $x$ below $x^{\ell+1}$ in its expansion at the boundary.
In \cite[Thm. 4.5]{My Preprint} we prove that
$\FP_{\varepsilon =0} \int_{x=\varepsilon} \Fun =0$ under the (slightly) more general assumption that the  expansion of $h_x$ is even in $x$ below $x^{n-1}$ and the $x^{n-1}$ term in its expansion is trace-free. This is true for example on Poincar\'e-Einstein manifolds, of particular interest since they occur in the AdS/CFT correspondence in physics.
(Note that in \cite{My Preprint} we state this theorem for asymptotically hyperbolic manifolds, i.e. those with $\alpha \equiv 1$, but the same proof works for conformally compact manifolds.) We state this formally as a corollary to Theorem \ref{ChernRGB}.

\begin{corollary}
For conformally compact metrics that are even below $x^m$ for any special bdf, $x$,
\begin{equation*}
	{}^R\int \mathrm{Pff} = \chi(M).
\end{equation*}
In particular, from Corollary \ref{ConfComp1}, this implies that
\begin{equation*}
	(-1)^{m/2}  \lrpar{ {}^R\dim \mathcal{H}_{L^2}^{m/2}}
	= \sum_{k \geq \frac{m}{2}} (-1)^k \dim H^k(M,\partial M)
	 - \sum_{k<\frac{m}{2}} (-1)^k \dim H^k(M,\partial M).
\end{equation*}
\end{corollary}

An important invariant in the physical AdS/CFT theory is the renormalized volume,
\begin{equation*}
	{}^R\Vol = {}^R\int \dvol.
\end{equation*}
This depends on the choice of bdf used to renormalize the integral, but gives the same answer for every choice of special bdf. This follows from the fact that the expansion of $\dvol$ in $x$ consists of even terms up to $x^n$. In \cite{My Preprint} we show that the same is true for any scalar Riemannian invariant. This includes all of the heat invariants, i.e. the coefficients occurring in the short-time asymptotic expansion of the trace of the heat kernel.  In this section we will show that if the metric is even enough, the trace of the heat kernel for any fixed time $t>0$ itself has a renormalized integral independent of the choice of special bdf used to renormalize it.

To this end, 
we will use the representation of the heat kernel as the inverse Laplace transform of the resolvent. Then mapping properties of the resolvent will be parlayed into information about its expansion in $x$, a special bdf, culminating with the evenness of the heat kernel at the front face up to $\rho_{11}^m$.

%%%%%%%%%%%%%%%%%%%%%%%%%%%%%%%%%%%%%%%%%%%%%%%%%%%%%%%%%%%%%%%%%%%%%%%%%%%%%%%%%%%%%%%%%%%%%%%%%%%%%%%%%%%%%%%%%%%%%%%%%%%%
\subsection{Even Functions and Operators on the Stretched Double Space} \label{EvenSubCalc} $ $\newline
%%%%%%%%%%%%%%%%%%%%%%%%%%%%%%%%%%%%%%%%%%%%%%%%%%%%%%%%%%%%%%%%%%%%%%%%%%%%%%%%%%%%%%%%%%%%%%%%%%%%%%%%%%%%%%%%%%%%%%%%%%%%

We proceed as in \cite[Chapter 7]{APS Book} to define even pseudodifferential operators.
We start by recalling, as in \cite[$\mathcal{x}$2]{Guillarmou} the space of even functions on an asymptotically hyperbolic manifold.
We use this to define a space of even functions and then operators on the double stretched product. Our goal is eventually to show that the heat kernel is in this calculus (cf. Corollary \ref{HeatKerFun}).

Given a special bdf, $x$, we can use the flow of the gradient $\nabla_{x^2g}x$ to identify a neighborhood of the boundary, $\curly{U}_x$ with a product neighborhood, $[0,\varepsilon_x)\times \partial M$. In the interior of this neighborhood, $g$ can be expressed by
\begin{equation}\label{MetricNorForm}
	\frac{dx^2}{x^2} + \frac{h(x,y,dy)}{x^2}.
\end{equation}
We will assume that the expansion of $h$ at the boundary contains only even powers of $x$ below $x^{2\ell}$, and say that $g$ is even mod $x^{2\ell}$.
It turns out \cite[Lemma 2.1]{Guillarmou} that the coordinate changes $(x,y) \mapsto (\tilde{x}, \tilde{y})$ that preserve the form of the metric on $[0,\varepsilon_x)\times \partial M$ have local expansions at the boundary of the form
\begin{equation}\label{CoordChange}
	\tilde{x} = x\sum_{j=0}^{\ell+1} a_j(y)x^{2j} + \curly{O}\lrpar{x^{2\ell+4}}, \phantom{xx}
	\tilde{y} = \sum_{j=0}^{\ell+1} b_j(y)x^{2j} + \curly{O}\lrpar{x^{2\ell+3}}.
\end{equation}
In particular, if $g$ is even mod $x^{2\ell}$ for one special bdf, it is even mod $\tilde{x}^{2\ell}$ for any special bdf.
We will refer to coordinates for a neighborhood of $\partial M$ of this type as `special' coordinates.

Similarly, the spaces of even functions and odd functions mod $x^{2\ell}$, respectively denoted $C^{\infty}_{even}(M)$ and $C^{\infty}_{odd}(M)$ are also well-defined independently of the choice of $x$.
We want to define $C^{\infty}_{even}(M^2_0)$. Whatever this space is, it should certainly contain 
\begin{equation}\label{subspace}
	\beta_L^*\lrpar{C^{\infty}_{even}(M)} \cdot \beta_{R}^*\lrpar{C^{\infty}_{even}(M)}
	+ \beta_L^*\lrpar{C^{\infty}_{odd}(M)} \cdot \beta_{R}^*\lrpar{C^{\infty}_{odd}(M)}.
\end{equation}

We start with the following polar coordinates on $M^2_0$
\begin{equation}\label{PolarCoords}
	\lrpar{ R, \omega, v } := 
		\lrpar{ \lrpar{x^2 + (x')^2 + |y-y'|^2}^{\frac12}, \frac1R\lrpar{x,\frac{y-y'}2,x'}, \frac{y+y'}2 },
\end{equation}
where $x$, $x'$ are the same special bdf on their respective factors. 
Consider, for $f \in C^{\infty}_{even}(M)$,
\begin{equation*}
	\beta_L^*\lrpar{f} = f\lrpar{R\omega_0, v+R\omega'} 
		\sim \sum_{even} \lrpar{R\omega_0}^k f_k'\lrpar{v+R\omega'}
		\sim \sum R^k f_k \lrpar{ \omega_0, \omega', v}. 
\end{equation*}
Because $f$ is $even \mod x^{2\ell}$ we can conclude that below $R^{2\ell}$, the even terms in this expansion are even with respect to the reflection
\begin{equation*}
	\omega' \xrightarrow{\Phi} -\omega'
\end{equation*}
while the odd terms are odd with respect to $\Phi$.
The same is true for
\begin{equation*}
	\beta_R^*\lrpar{f} = f\lrpar{R\omega_m, v-R\omega'}.
\end{equation*}

Thus we can achieve the inclusion of \eqref{subspace} by defining $F \in C^{\infty}_{even}\lrpar{M^2_0}$ if, in the coordinates \eqref{PolarCoords},
\begin{equation}\label{EvenStretch}
	F(R, \omega, v) \sim \sum_{j < 2\ell} R^j F_j(\omega, v) + R^{2\ell} F'(R, \omega, v)
\end{equation}
with $F_j$ respectively even or odd with respect to $\Phi$ if $j$ is even or odd. Similarly, $F \in C^{\infty}_{odd}\lrpar{M^2_0}$ if $F_j$ is even with respect to $\Phi$ for odd $j$, and odd with respect to $\Phi$ for even $j$.
Notice that functions that are odd with respect to $\Phi$ necessarily vanish at the diagonal. Hence, if $F$ is an even function on $M^2_0$ and we identify the diagonal with $M$, $F$ restricts to the diagonal to an element of $C^{\infty}_{even}(M)$.

For this space to be well-defined, it should be independent of the coordinates on $M$ that we started with.
So consider $(\tilde{R}, \tilde{\omega}, \tilde{v})$ defined by \eqref{PolarCoords} in terms of $(\tilde{x}, \tilde{y})$ satisfying \eqref{CoordChange}. 
As in \cite[Proposition 7.7]{APS Book}, note that
\begin{equation*}\begin{split}
	\lrpar{\frac{\tilde{R}}{R}}^2 &=\frac{\tilde{x}^2 + \tilde{x'}^2+ |\tilde{y}-\tilde{y'}|^2}{R^2}\\
	&=\omega_0^2 \lrpar{\sum a_j\lrpar{v + R\omega'}\lrpar{R\omega_0}^{2j} }^2
	+\omega_m^2 \lrpar{\sum a_j\lrpar{v - R\omega'}\lrpar{R\omega_m}^{2j} }^2
	+ \left\lvert \omega' \cdot G\lrpar{y,y'} \right\rvert^2
\end{split}\end{equation*}
with $G$ a smooth matrix with non-zero determinant for $y$ near $y'$, hence $\tilde{R}\in C^\infty_{odd}\lrpar{M^2_0}$. Also since $\tilde{x}$, $\tilde{y}$ are even and odd functions on $M$ respectively, they lift to even and odd functions on $M^2_0$, together with $\tilde{R}$ odd, this implies
\begin{equation*}
	\tilde{\omega}_0, \tilde{\omega}_m, \tilde{v} \in C^\infty_{even}\lrpar{M^2_0}, \phantom{xxx}
	\tilde{R}, \tilde{\omega}' \in C^\infty_{odd}\lrpar{M^2_0},
\end{equation*}
which in turn shows that we get the same sets of even and odd functions on $M^2_0$ starting with $(x,y)$ or with $(\tilde{x}, \tilde{y})$, i.e. these sets are defined independently of choice of coordinates (from among special coordinates).

For future reference, we note that we have a similar expansion for any other bdf for the front face $\rho \in C^{\infty}_{odd}(M)$. Thus, for instance given coordinates $(x,y)$ as above, we can use projective coordinates on $M^2_0$ away from $\{ x =0 \}$ of the form
\begin{equation}\label{s'u'}
	\lrpar{ x,y,s',u' } := \lrpar{ x,y, \frac{x'}x, \frac{y'-y}x }
\end{equation}
and a function $f \in C^{\infty}_{even}(M^2_0)$ will have an expansion
\begin{equation*}
	f \sim \sum_{j < 2\ell} x^j f_j\lrpar{y,s',u'} + x^{2\ell}f'\lrpar{x,y,s',u'}
\end{equation*}
with $f_j$ respectively even or odd with respect to $u' \mapsto -u'$ for $j$ even or odd.

We next extend the definition of even functions to even operators. 
If $u$ is a distribution on $M^2_0$, we let $\Phi$ act on $u$ by demanding that for any test function, $\phi$,
\begin{equation*}
	\left< \Phi^*u, \phi \right> = \left< u, \phi \circ \Phi \right>.
\end{equation*}
Now consider the expansion of
an element $K \in \curly{A}^\curly{E}I^k\lrpar{M^2_0, \Omega^{1/2}}$ at the front face
\begin{equation}\label{Exp11}
	\curly{K} \sim \sum_{\lrpar{z,p}\in E_{11}} \sum_{q \leq p} 
		\rho_{11}^z \lrpar{\log \rho_{11}}^{\ell} K_{\lrpar{z,q}},
\end{equation}
where, with $\curly{E}' = \lrpar{ E_{01}, E_{10}}$,
\begin{equation*}
	K_{\lrpar{z,q}} \in \curly{C}_c^{\infty}
		\lrpar{ [0,1); \curly{A}_{phg}^{\curly{E}'}I^{k+1/4} %
		\lrpar{\bhs{11}, \mathrm{diag}_0 \cap \bhs{11}}}.
\end{equation*}
We will say that $K$ is even if $E_{11} \cap \{\Re(z)<2\ell\} \subset \mathbb{N}_0$ and $\Phi^*K_{(j,0)} = \lrpar{-1}^jK_{(j,0)}$
We will denote the space of operators whose kernels are even at the front face in this sense by $\Psi^{k,\curly{E}}_{0,even}(M)$ and the corresponding space of odd operators by $\Psi^{k,\curly{E}}_{0,odd}(M)$.

If $K$ is given by a function (not just a distribution) on $M^2_0$, then the criteria for the even (resp. odd) subcalculus is the same as on smooth functions, \eqref{EvenStretch}.
On the other hand, if $A$ is a $0$-differential operator then it has the expected behavior on the $even$ operators.

\begin{lemma}\label{DifferentialComp}
The operators $x\partial_x$ and $x\partial_{y_i}$ are respectively $even$ and $odd$ operators.
Composition with $x\partial_x$ preserves even and odd operators, while composition with $x\partial_{y_i}$ interchanges even and odd operators.
\end{lemma}

\begin{proof}
We will work in the coordinates
\begin{equation*}
	\lrpar{ s, u, x', y'} := \lrpar{ \frac{x}{x'}, \frac{y-y'}{x'}, x', y' }
\end{equation*}
on $M^2_0$. Assume $\curly{K}$ is a distribution on $M^2_0$ such that
$\Phi_*\curly{K} = \lrpar{-1}^j \curly{K}$.

For $x\partial_x$ we find
\begin{equation*}\begin{split}
	\left< \Phi_*\lrpar{s\partial_s\curly{K}}, \beta^*\lrpar{ \phi \otimes \psi} \right>
	&= - \int  \curly{K} \partial_s \phi\lrpar{ sx', -x'u+y'} \psi \lrpar{x',y'}
		ds du \frac{dx' dy'}{(x')^m} \\
	&= - \int  \curly{K} x'\partial_x'\phi\lrpar{ sx', -x'u+y'} \psi \lrpar{x',y'}
		ds du \frac{dx' dy'}{(x')^m} \\
	&= - \int  \curly{K} \Phi^*x'\partial_x'\phi\lrpar{ sx', x'u+y'} \psi \lrpar{x',y'}
		ds du \frac{dx' dy'}{(x')^m} \\
	&= -(-1)^j \int  \curly{K} x'\partial_x'\phi\lrpar{ sx', x'u+y'} \psi \lrpar{x',y'}
		ds du \frac{dx' dy'}{(x')^m} \\
	&= (-1)^j \left< \lrpar{s\partial_s\curly{K}}, \beta^*\lrpar{ \phi \otimes \psi} \right>
\end{split}\end{equation*}
and the same computation in the coordinates \eqref{s'u'} shows that
\begin{equation*}
	\left< \Phi_*\lrpar{\curly{K}\circ s'\partial_{s'}}, \beta^*\lrpar{ \phi \otimes \psi} \right>
	= (-1)^j \left< \curly{K}\circ s'\partial_{s'}, \beta^*\lrpar{ \phi \otimes \psi} \right>.
\end{equation*}
Similarly, for $x\partial_{y_i}$,
\begin{equation*}\begin{split}
	\left< \Phi_*\lrpar{s\partial_{u_i}\curly{K}}, \beta^*\lrpar{ \phi \otimes \psi} \right>
	&= - \int  \curly{K} \partial_{u_i} \phi\lrpar{ sx', -x'u+y'} \psi \lrpar{x',y'}
		ds du \frac{dx' dy'}{(x')^m} \\
	&= - \int  \curly{K} (-x')\partial_{y_i'}\phi\lrpar{ sx', -x'u+y'} \psi \lrpar{x',y'}
		ds du \frac{dx' dy'}{(x')^m} \\
	&=  \int  \curly{K} \Phi^*x'\partial_{y_i'}\phi\lrpar{ sx', x'u+y'} \psi \lrpar{x',y'}
		ds du \frac{dx' dy'}{(x')^m} \\
	&= (-1)^j \int  \curly{K} x'\partial_{y_i'}\phi\lrpar{ sx', x'u+y'} \psi \lrpar{x',y'}
		ds du \frac{dx' dy'}{(x')^m} \\
	&= (-1)^{j+1} \left< \lrpar{s\partial_{u_i}\curly{K}}, \beta^*\lrpar{ \phi \otimes \psi} \right>
\end{split}\end{equation*}
and likewise
\begin{equation*}
	\left< \Phi_*\lrpar{\curly{K}\circ s'\partial_{u_i'}}, \beta^*\lrpar{ \phi \otimes \psi} \right>
	= (-1)^{j+1} \left< \curly{K}\circ s'\partial_{u_i'}, \beta^*\lrpar{ \phi \otimes \psi} \right>.
\end{equation*}

From these computations, with $\curly{K}$ the kernel of the identity we see that 
\begin{equation*}
	x\partial_x \in {}^0\Psi^1_{even}\lrpar{M;\Omega^{\frac12}}, \phantom{x}
	x\partial_x \in {}^0\Psi^1_{odd}\lrpar{M;\Omega^{\frac12}}.
\end{equation*}

Taking instead $\curly{K}$ as the $j^{\mathrm{th}}$ term in the expansion of an operator at the front face, we see that $x\partial_x$ preserves the parity with respect to $\Phi$ while $x\partial_{y_i}$ reverses it, proving the rest of the lemma.
\end{proof}

We prove a composition result, first for kernels in $\Psi^{-\infty}_{0,even}$.
\begin{lemma} \label{ResidualComp}
\begin{equation*}
	\Psi^{-\infty}_{0,even}(M) \circ \Psi^{-\infty}_{0,even}(M) \subset \Psi^{-\infty}_{0,even}(M)
\end{equation*}
\end{lemma}
\begin{proof}
Let $A,B \in \Psi^{-\infty}_{0,even}(M)$, recall that
\begin{equation*}
	\curly{K}_{A \circ B} = \lrpar{\beta_{LR}}_*
		\lrpar{\beta_{LM}^*\lrpar{\curly{K}_A}\beta_{MR}^*\lrpar{\curly{K}_B}}.
\end{equation*}
Let $(x, y)$ be special coordinates on the left factor with identical coordinates $(\hat{x},\hat{y})$ and $(x',y')$ on the middle and right factors.
We introduce projective coordinates
\begin{equation*}
	\lrpar{x,y,s',u'} := \lrpar{x,y,\frac{x'}x,\frac{y'-y}x}, \phantom{x}
	\lrpar{x,y,\hat{s},\hat{u}} := \lrpar{x,y,\frac{\hat{x}}x,\frac{\hat{y}-y}x}, \phantom{x}
	\lrpar{s'x, y+u'x, \frac{\hat{s}}{s'}, \frac{\hat{u}-u'}{s'}}
\end{equation*}
on the double spaces $LR$, $LM$, and $MR$ respectively. In these coordinates, neglecting a factor of
\begin{equation*}
			\left\lvert \frac{dxdy}{x^m}\frac{ds'du'}{\lrpar{s'}^m} \right\rvert^{\frac12}
\end{equation*}
on both sides,
\begin{equation*}
	\curly{K}_{A \circ B}\lrpar{x,y,s',u'} 
	=\int
	\curly{K}_{A}\lrpar{x,y,\hat{s},\hat{u}} 
	\curly{K}_{B}\lrpar{s'x, y+u'x, \frac{\hat{s}}{s'}, \frac{\hat{u}-u'}{s'}}
		\frac{d\hat{s}d\hat{u}}{\lrpar{s'}^m}.
\end{equation*}
We know that $A$ is an $even$ operator, hence
\begin{equation}\label{AExp}
	\curly{K}_{A}\lrpar{x,y,\hat{s},\hat{u}} 
	\sim \sum_j x^j \lrpar{\curly{K}_A}_j\lrpar{y,\hat{s},\hat{u}} + \curly{O}\lrpar{x^{2\ell}}
\end{equation}
with $\Phi^*_{\hat{u}}\lrpar{\curly{K}_A}_j = \lrpar{-1}^j\lrpar{\curly{K}_A}_j$, and similarly for $\curly{K}_B$.

We claim that 
\begin{equation}\label{BClaim}
	\curly{K}_{B}\lrpar{s'x, y+u'x, \frac{\hat{s}}{s'}, \frac{\hat{u}-u'}{s'}}
	\sim \sum_k x^k \lrpar{\curly{L}_B}_k\lrpar{y,\hat{s},s',u',\hat{u}} + \curly{O}\lrpar{x^{2\ell}}
\end{equation}
with $\Phi^*_{u'}\lrpar{\curly{L}_B}_j = \lrpar{-1}^j\Phi^*_{\hat{u}}\lrpar{\curly{L}_B}_j$.
Indeed, the analogue of \eqref{AExp} for $B$ yields
\begin{equation*}\begin{split}
	\curly{K}_{B}\lrpar{s'x, y+u'x, \frac{\hat{s}}{s'}, \frac{\hat{u}-u'}{s'}}
	&\sim \sum_j (s'x)^j \lrpar{\curly{K}_B}_j\lrpar{y+u'x, \frac{\hat{s}}{s'}, \frac{\hat{u}-u'}{s'}} 
	+ \curly{O}\lrpar{x^{2\ell}}\\
	&=: \sum_j (s'x)^j \lrspar{%
	\sum_i x^i \curly{K}_{j,i}\lrpar{y,u', \frac{\hat{s}}{s'}, \frac{\hat{u}-u'}{s'}} }
	+ \curly{O}\lrpar{x^{2\ell}},
\end{split}\end{equation*}
with $K_{j,i}$ respectively even or odd with respect to the first $u'$ in its arguments when $i$ is even or odd. If we write 
\begin{equation*}
	 \sum_j (s'x)^j \lrspar{%
	\sum_i x^i \curly{K}_{j,i}\lrpar{y,u', \frac{\hat{s}}{s'}, \frac{\hat{u}-u'}{s'}} }
	= \sum_j \lrpar{s'x}^j \lrspar{%
	\sum_i x^i K_{j,i}\lrpar{y,s',\hat{s},u',\hat{u} } },
\end{equation*}
then
\begin{equation*}
	K_{j,i}\lrpar{y,s',\hat{s},-u',-\hat{u} } 
	=\curly{K}_{j,i}\lrpar{y,-u', \frac{\hat{s}}{s'}, -\frac{\hat{u}-u'}{s'}}	
	=\lrpar{-1}^{j+i} K_{j,i}\lrpar{y,s',\hat{s},u',\hat{u} }. 
\end{equation*}
So with 
\begin{equation*}
	\lrpar{\curly{L}_B}_k\lrpar{y,\hat{s},s',u',\hat{u}} =  
	\sum_{i+j=k} K_{j,i}\lrpar{y,s',\hat{s},u',\hat{u} } 
\end{equation*}
we have
\begin{equation*}
	\lrpar{\curly{L}_B}_k\lrpar{y,\hat{s},s',-u',\hat{u}} =  
	(-1)^k \lrpar{\curly{L}_B}_k\lrpar{y,\hat{s},s',u',-\hat{u}},
\end{equation*}
hence the claim.

Finally, we can finish the proof of the lemma by noting
\begin{equation*}\begin{split}
	\curly{K}_{A \circ B}\lrpar{x,y,s',u'} 
	&=\int
	\curly{K}_{A}\lrpar{x,y,\hat{s},\hat{u}} 
	\curly{K}_{B}\lrpar{s'x, y+u'x, \frac{\hat{s}}{s'}, \frac{\hat{u}-u'}{s'}}
		\frac{d\hat{s}d\hat{u}}{\lrpar{s'}^m} \\
	&\sim \sum_{j,k} x^{j+k} \int
	\lrpar{\curly{K}_{A}}_j\lrpar{y,\hat{s},\hat{u}} 
	\lrpar{\curly{L}_{B}}_k\lrpar{y,\hat{s},s',u',\hat{u}}
		\frac{d\hat{s}d\hat{u}}{\lrpar{s'}^m} + \curly{O}\lrpar{x^{2\ell}}
\end{split}\end{equation*}
and
\begin{equation*}\begin{split}
	\Phi^*_{u'} 
	\lrpar{\curly{K}_{A\circ B}}_k\lrpar{y,s',u'}
	&= \sum_{i+j=k} \int
	\lrpar{\curly{K}_{A}}_i\lrpar{y,\hat{s},\hat{u}} 
	\lrpar{\curly{L}_{B}}_j\lrpar{y,\hat{s},s',-u',\hat{u}}
		\frac{d\hat{s}d\hat{u}}{\lrpar{s'}^m} \\
	&= \sum_{i+j=k} (-1)^j\int
	\lrpar{\curly{K}_{A}}_i\lrpar{y,\hat{s},\hat{u}} 
	\lrpar{\curly{L}_{B}}_j\lrpar{y,\hat{s},s',u',-\hat{u}}
		\frac{d\hat{s}d\hat{u}}{\lrpar{s'}^m}\\
	&= \sum_{i+j=k} (-1)^j\int
	\lrpar{\curly{K}_{A}}_j\lrpar{y,\hat{s},-\hat{u}} 
	\lrpar{\curly{L}_{B}}_k\lrpar{y,\hat{s},s',u',\hat{u}}
		\frac{d\hat{s}d\hat{u}}{\lrpar{s'}^m}\\
	&=\sum_{i+j=k} (-1)^{i+j}\int
	\lrpar{\curly{K}_{A}}_i\lrpar{y,\hat{s},\hat{u}} 
	\lrpar{\curly{L}_{B}}_j\lrpar{y,\hat{s},s',u',\hat{u}}
		\frac{d\hat{s}d\hat{u}}{\lrpar{s'}^m} \\
	&=\lrpar{-1}^{k} \lrpar{\curly{K}_{A\circ B}}_k\lrpar{y,s',u'}.
\end{split}\end{equation*}
\end{proof}

It is now easy to extend composition to certain distributional kernels. Recall that if the composition is defined (i.e., $\mathrm{Re}\lrpar{E_{01}} + \mathrm{Re}\lrpar{F_{10}} > -1$)
\begin{equation*}
	{}^0\Psi^{j,\curly{E}}_{\mathrm{even}}(M;\hat{\rho}) \circ
	{}^0\Psi^{k,\curly{F}}_{\mathrm{even}}(M;\hat{\rho}) \subset
	{}^0\Psi^{j+k,\curly{G}}(M)
\end{equation*}
where
\begin{equation}\label{IndSets}\begin{split}
	G_{10}&=\lrpar{E_{11} + F_{10}} \bar{\cup} E_{10}, \phantom{xxx}
	G_{01}=\lrpar{E_{01} + F_{11}} \bar{\cup} F_{01}, \\
	G_{11}&= \lrpar{E_{11} + F_{11}} \bar{\cup} \lrpar{E_{10} + F_{01} + m}.
\end{split}\end{equation}

\begin{proposition}\label{Composition}
If $A \in \Psi^{j,\curly{E}}_{0, even}(M)$, $B \in \Psi^{k,\curly{F}}_{0, even}(M)$ satisfy
\begin{itemize}
	\item $\mathrm{Re}\lrpar{E_{01}} + \mathrm{Re}\lrpar{F_{10}} > -1$
	\item $\Re\lrpar{ E_{10}+F_{01}+m } > 2\ell$,
\end{itemize}
then $A \circ B \in \Psi^{j+k,\curly{G}}_{0, even}(M)$.
\end{proposition}

\begin{proof}
Because $\Re\lrpar{ E_{10}+F_{01}+m } > 2\ell$, the expansion of $\curly{K}_{A\circ B}$ at the front face below order $2\ell$ comes from the expansions of $\curly{K}_A$, $\curly{K}_B$ at the front face, so it suffices to prove that even elements of the small calculus compose.

If $A$ and $B$ are in the small calculus and $j$, $k$ are sufficiently negative then their kernels restrict to the diagonal as functions and the proof of composition in Lemma \ref{ResidualComp} extends verbatim to cover this case.

If the orders of $A$ or $B$ are not sufficiently negative to apply Lemma \ref{ResidualComp} then we use Lemma \ref{DifferentialComp} to reduce the order of the operators involved.
%Any pseudodifferential operator can be written as a differential operator applied to a lower order pseudodifferential operator.
For instance, if $C \in {}^0\Psi^k_{even}$ is any $even$ $0$-differential operator which is invertible as a $0$-pseudodifferential operator (e.g. $\Delta +\mathrm{Id}$ which we will soon see is $even$).
Since every $0$-pseudodifferential operator can be written as an $even$ operator plus an $odd$ operator and we know that composition with $C$ preserves $even$ and $odd$ operators, it follows that $C^{-1}$ must also preserve $even$ and $odd$ operators in the small calculus (even if we do not know that $C^{-1}$ is even). Hence writing $A = C^k\lrpar{C^{-k}A}$ and similarly for $B$ we can write the composition of $A$ and $B$ as a composition of lower-order $even$ pseudodifferential operators with $even$ differential operators. Thereby reducing the problem to Lemmas  \ref{DifferentialComp} and \ref{ResidualComp}.
\end{proof}

We can also check that $even$ operators preserve $even$ functions. For instance, if $K_A$ is an even operator kernel in the small calculus and $f$ an even function on $M$, then in the coordinates \eqref{s'u'} we can split $K=K_e+xK_o$ mod $\curly{O}(x^{2\ell})$ with $K_e$, $K_o$ even in $x$ and respectively even and odd with respect to $\Phi$, and similarly $\beta^*f = f_e + xf_o$. Then we can see that $A(f)$ is an even function,
\begin{equation}\begin{split}\label{EvenChar}
	\int K\lrpar{x,y,s',u'}&f(xs', xu'+y) \frac{ds'du'}{(s')^m}\\
	&= \int \lrpar{K_{e}+xK_{o}}(x^2,y,s',u')\lrpar{f_{e}+xf_{o}}(x^2,y,s',u')
	 + \curly{O}\lrpar{x^{2\ell}} \frac{ds'du'}{(s')^m} \\
	&= \int \lrpar{K_{e}f_{e} + x^2K_{o}f_{o}}(x^2,y,s',u') \frac{ds'du'}{(s')^m}
	 +  \curly{O}\lrpar{x^{2\ell}},
\end{split} \end{equation}
where we have used that the integral with respect to $u'$ of an odd function in $u'$ vanishes.
Similarly, if $K$ is odd and $f$ even then $K(f)$ is odd.

Finally, we point out that the even subcalculus has distinguished image under the symbol maps. 
The symbol of a $0$-pseudodifferential operator is a section of the zero cotangent bundle and the symbol of an $even$ $0$-pseudodifferential operator will be a $even$ section. 
Similarly, after we identify the normal operator at the fiber over a point $q \in \partial M$ with an operator on 
$\mathbb{R}^+ \times \mathbb{R}^{m-1}$, it is clear that
the image of the even $0$-pseudodifferential operators on $M$ consists of the even $0$-pseudodifferential operators on $\mathbb{R}^+ \times \mathbb{R}^{m-1}$.
Note that these maps are surjective and their kernels consist of $even$ $0$-pseudodifferential operators.

%%%%%%%%%%%%%%%%%%%%%%%%%%%%%%%%%%%%%%%%%%%%%%%%%%%%%%%%%%%%%%%%%%%%%%%%%%%%%%%%%%%%%%%%%%%%%%%%%%%%%%%%%%%%%%%%%%%%%%%%%%%%
\subsection{The Laplacian on forms} $ $\newline
%%%%%%%%%%%%%%%%%%%%%%%%%%%%%%%%%%%%%%%%%%%%%%%%%%%%%%%%%%%%%%%%%%%%%%%%%%%%%%%%%%%%%%%%%%%%%%%%%%%%%%%%%%%%%%%%%%%%%%%%%%%%

In this section we will identify $even$ forms on $M$ and $even$ operators between them. We will see that the Hodge Laplacian is such an operator and in the next section show that its resolvent and heat kernel are as well.
For clarity, we start by identifying even one-forms in ${}^0\Lambda^1\lrpar{M}$,
\begin{equation*}
	\omega = b(x,y) \frac{dx}x + \sum_{i<m} a_i(x,y) \frac{dy^i}{x}.
\end{equation*}
Since $\frac{dx}x$ and $\frac{dy}x$ are even and odd respectively with regards to (formally) replacing $x$ with $-x$,
we will say that $\omega$ is an even form if $b \in C^{\infty}_{even}(M)$ and $a \in C^{\infty}_{odd}(M)$.
A moment's thought shows the consistency of this definition in that
\begin{equation*}
	C^{\infty}_{even}\lrpar{M;\Omega^{\frac12}} \xrightarrow{{}^0d} 
		{}^0\Lambda^1_{even}\lrpar{M;\Omega^{\frac12}}.
\end{equation*}

Similarly, a section of ${}^0\Lambda^k\lrpar{M}$ can be written in local coordinates as
\begin{equation*}
	\omega = \sum_{|\alpha|=k} a_{\alpha}(x,y) \frac{dy^{\alpha}}{x^k}
		+\sum_{|\beta|=k-1} b_{\beta}(x,y) \frac{dx}x\wedge\frac{dy^{\beta}}{x^{k-1}},
\end{equation*}
where $\alpha$ and $\beta$ are multi-index sets from $\{1,\ldots, m-1 \}$.
Just as for one-forms, we will say that $\omega$ is an $even$ form, 
\begin{equation}\label{EvenForm}
	\omega \in {}^0\Lambda^k_{even}\lrpar{M}, 
\end{equation}
if $\alpha \in C^{\infty}_{even}(M)$ and $\beta \in C^{\infty}_{odd}(M)$ for even $k$, while for odd $k$, we will demand
$\alpha \in C^{\infty}_{odd}(M)$ and $\beta \in C^{\infty}_{even}(M)$. Forms in ${}^0\Lambda^k_{odd}\lrpar{M}$ are defined analogously so that
\begin{align}
	{}^0\Lambda^k(M) = {}^0\Lambda^k_{even}(M) + {}^0\Lambda^k_{odd}(M) \label{LambdaSum} \\
	{}^0\Lambda^k_{even}(M) \cap {}^0\Lambda^k_{odd}(M) 
		= x^{2\ell} \phantom{x} {}^0\Lambda^k(M).
\end{align}

Notice that if $(\tilde{x}, \tilde{y})$ are another set of special coordinates as in \eqref{CoordChange}, then
\begin{equation}\label{dxChange}
	\frac{d\tilde{x}}{\tilde{x}} = \frac{1}{\tilde{x}}
	\lrpar{\frac{\partial \tilde{x}}{\partial x} dx 
		+ \sum_{i<m} \frac{\partial \tilde{x}}{\partial {y_i}} dy^{i}}
		= b(x,y) \frac{dx}x + \sum_{i<m} a_i(x,y) \frac{dy^i}{x},
\end{equation}
with $b \in C^{\infty}_{even}(M)$ and $a_i \in C^{\infty}_{odd}(M)$ and similarly,
\begin{equation}\label{dyChange}
	\frac{d\tilde{y}^i}{\tilde{x}} 
		= b'(x,y) \frac{dx}x + \sum_{i<m} a'_i(x,y) \frac{dy^i}{x},
\end{equation}
with $b' \in C^{\infty}_{odd}(M)$ and $a'_i \in C^{\infty}_{even}(M)$. Hence the spaces of $odd$ and $even$ forms as defined here are independent of the choice of special coordinates. Note in contrast that the splitting of the bundle of $0$-differential forms into tangential and normal parts is only well-defined up to first order in $x$.

The kernel of a $0$-pseudodifferential operator from $j$-forms to $k$-forms is an element of
\begin{equation*}
	{}^0\Psi^{s,\curly{E}}\lrpar{M;{}^0\Omega^{1/2}} \otimes_{C^{\infty}\lrpar{M^2_0}}
	C^{\infty}\lrpar{M^2_0;%
	\beta^*\Hom\lrpar{{}^0\Lambda^j_L \otimes {}^0\Omega^{-1/2}, %
		{}^0\Lambda^k_R \otimes {}^0\Omega^{-1/2} }}.
\end{equation*}
Neglecting density terms, we can write any such operator as 
\begin{equation}\label{HomSum} \begin{split}
	F &=
	\sum_{\alpha,\beta} f^{nn}_{\alpha,\beta}\lrpar{R,\omega,v}
		\beta_{L}^*\lrpar{\frac{dx}x\frac{dy^{\alpha}}x}\beta_{R}^*\lrpar{\frac{dx}x\frac{dy^{\beta}}x}'
	+\sum_{\alpha,\beta} f^{tn}_{\alpha,\beta}\lrpar{R,\omega,v}
		\beta_{L}^*\lrpar{\frac{dy^{\alpha}}x}\beta_{R}^*\lrpar{\frac{dx}x\frac{dy^{\beta}}x}' \\
	&\phantom{x} +\sum_{\alpha,\beta} f^{nt}_{\alpha,\beta}\lrpar{R,\omega,v}
		\beta_{L}^*\lrpar{\frac{dx}x\frac{dy^{\alpha}}x}\beta_{R}^*\lrpar{\frac{dy^{\beta}}x}'
	+\sum_{\alpha,\beta} f^{tt}_{\alpha,\beta}\lrpar{R,\omega,v}
		\beta_{L}^*\lrpar{\frac{dy^{\alpha}}x}\beta_{R}^*\lrpar{\frac{dy^{\beta}}x}',
\end{split}\end{equation}
with $f^{tt}_{\alpha,\beta}$, $f^{tn}_{\alpha,\beta}$, $f^{nt}_{\alpha,\beta}$, $f^{nn}_{\alpha,\beta}$ elements of 
${}^0\Psi^{s,\curly{E}}\lrpar{M;{}^0\Omega^{1/2}}$
Define 
\begin{equation*}
	{}^0\Psi^{s,\curly{E}}_{even}\lrpar{M;{}^0\Lambda^j,{}^0\Lambda^k}
\end{equation*}
as those operators with 
\begin{equation*}
	f^{nn}_{\alpha,\beta}, f^{tt}_{\alpha,\beta} 
		\in {}^0\Psi^{s,\curly{E}}_{even}\lrpar{M;{}^0\Omega^{1/2}}, \phantom{xx}
	f^{nt}_{\alpha,\beta}, f^{tn}_{\alpha,\beta} 
		\in {}^0\Psi^{s,\curly{E}}_{odd}\lrpar{M;{}^0\Omega^{1/2}}
\end{equation*}
for any $\alpha$ and $\beta$. Note that operators of this form compose with suitable restrictions on the index sets $\curly{E}$. 

We need to check coordinate invariance of this space. It suffices to show that for $(\tilde{x}, \tilde{y})$
any other special coordinates, $\beta_L^*\lrpar{\frac{d\tilde{x}}{\tilde{x}}}$ and  $\beta_L^*\lrpar{\frac{d\tilde{y}_i}{\tilde{x}}}$ are respectively even and odd with respect to $\Psi$. This follows from \eqref{dxChange}, \eqref{dyChange} and the fact that 
$\beta^*_LC^{\infty}_{even} \subset C^{\infty}_{even}(M^2_0)$, 
$\beta^*_LC^{\infty}_{odd} \subset C^{\infty}_{odd}(M^2_0)$. In the same way, we can verify the comforting fact that
\begin{equation*}
	 \lrpar{{}^0\Lambda^k_{even}}_L' \otimes \lrpar{{}^0\Lambda^k_{even}}_R
	+ \lrpar{{}^0\Lambda^k_{odd}}_L' \otimes \lrpar{{}^0\Lambda^k_{odd}}_R
\end{equation*}
is a subset of 
${}^0\Psi^{-\infty}_{even}\lrpar{M;{}^0\Lambda^k,{}^0\Lambda^k}$. 

\begin{lemma} \label{LapIsEven}
For any $0\leq j \leq m$,
\begin{itemize}
\item[a)] ${}^0d \in {}^0\Psi^{1}_{odd}\lrpar{M;{}^0\Lambda^j,{}^0\Lambda^{j+1}}$
\item[b)] $* \in {}^0\Psi^{0}_{odd}\lrpar{M;{}^0\Lambda^j,{}^0\Lambda^{m-j}}$
\item[c)] $\Delta \in {}^0\Psi^{2}_{even}\lrpar{M;{}^0\Lambda^j,{}^0\Lambda^j}$
\end{itemize}
\end{lemma}

\begin{proof}
{\em a)} With respect to the splitting of ${}^0\Lambda^j$ and ${}^0\Lambda^{j+1}$ into tangential and normal parts, the exterior derivative acts by
\begin{equation*}
	\begin{pmatrix}
	{}^0d_{\partial M} & 0 \\ x\partial_x & -{}^0d_{\partial M}
	\end{pmatrix},
\end{equation*}
and we have seen that $x\partial_x$, $x\partial_{y_i}$ are respectively even and odd $0$-pseudodifferential operators. \newline

{\em b)} With respect to the splitting of ${}^0\Lambda^j$ and ${}^0\Lambda^{m-j}$ into tangential and normal parts, the Hodge star acts by
\begin{equation*}
	\begin{pmatrix}
	0 & *_{\partial M}  \\ (-1)^j *_{\partial M} & 0
	\end{pmatrix},
\end{equation*}
and $*_{\partial M}$ acts by raising indices using the metric. Since the metric is even mod $x^{2\ell}$, its components lift to even functions on $M^2_0$, and thus 
$* \in C^{\infty}_{odd}\lrpar{M^2_0;\beta^*\Hom\lrpar{{}^0\Lambda^j_L ,{}^0\Lambda^{m-j}_R }}$. \newline

{\em c)} Follows from {\em (a)} and {\em (b)}.
\end{proof}

{\bf Remark.} Although the formulas \eqref{CoordChange} show that class of functions even mod $x^{2\ell+2}$ is invariantly defined, and similarly an `extended' even subcalculus replacing $x^{2\ell}$ by $x^{2\ell+2}$ in all of the constructions above, this does not contain the Laplacian. Indeed, it is easy to see (e.g., via \cite[(1.185)]{Besse}) that if we denote the $x^0$ and $x^{2\ell+1}$ terms in the metric by $\df{g}$ and $\df{h}$, the $x^{2\ell+1}$ term in the expansion of the Laplacian on functions (in a special bdf) is 
\begin{equation*}
	\df{g}\lrpar{\df{h}, \nabla_{\df{g}}d\cdot} - 
		\df{g}\lrpar{d\cdot, \mathrm{div}_{\df{g}}\df{h} +\frac12d\lrpar{\mathrm{tr}_{\df{g}}\df{h}}}.
\end{equation*}
For Poincar\'e-Einstein metrics this last term vanishes, nevertheless the first term contributes
\begin{equation*}
	x^{2\ell+1}h^{ij}x\partial_{y_i}x\partial_{y_j}
\end{equation*}
which is not in this `extended' even calculus!

%%%%%%%%%%%%%%%%%%%%%%%%%%%%%%%%%%%%%%%%%%%%%%%%%%%%%%%%%%%%%%%%%%%%%%%%%%%%%%%%%%%%%%%%%%%%%%%%%%%%%%%%%%%%%%%%%%%%%%%%%%%%
\subsection{The Resolvent is an even operator} \label{ResIsEven} $ $\newline
%%%%%%%%%%%%%%%%%%%%%%%%%%%%%%%%%%%%%%%%%%%%%%%%%%%%%%%%%%%%%%%%%%%%%%%%%%%%%%%%%%%%%%%%%%%%%%%%%%%%%%%%%%%%%%%%%%%%%%%%%%%%

In this section we will show that the construction of the resolvent yields an operator in the even subcalculus. There is something to check because composition in the even subcalculus is more restrictive than composition in the usual large calculus. We can eliminate this restriction by arranging for the index sets to vanish to high enough order at the side faces. Fortunately, this is controlled by the indicial roots and can be conveniently arranged.

\begin{lemma}\label{HypCase}
	On hyperbolic space $\mathbb{H}^m$, the Laplacian, resolvent and heat kernel are all in the even subcalculus.
\end{lemma}

\begin{proof}
It is possible to treat these operators very explicitly, see for instance \cite{CarronPedon}.
Because these operators are natural, they are invariant with respect to the isometries of the hyperbolic metric. These isometries are rich enough that their kernels, $\curly{K}(w,w')$ can only depend on the distance between the two points $w$, $w'$. It is known that hyperbolic distance, $\delta$, satisfies (e.g., \cite[6.6]{Mazzeo-Melrose})
\begin{equation*}
	\mathrm{cosh}\lrpar{\delta} = 1 + \frac{|w-w'|^2}{2xx'}
	= 1+ \frac{x^2 -2xx' + (x')^2 + |y-y'|^2}{2xx'}
	= \frac{x^2 + (x')^2 + |y-y'|^2}{2xx'}.
\end{equation*}
In `polar' coordinates \eqref{PolarCoords}
this pulls-back to the double space $\lrpar{\mathbb{H}^m}^2_0$ to
\begin{equation*}
	\beta^*\mathrm{cosh}\lrpar{\delta} = \frac{R^2}{2R^2\omega_0\omega_m}
	= \frac{1}{2\omega_0\omega_m}
\end{equation*}
as this is an $even$ function, the kernels will be in the even subcalculus. In fact, this shows that the kernels will have a constant expansion at the front face, and are globally invariant under $\Phi$.

\end{proof}

Recall (e.g., \cite[(4.15)]{Grieser}) that the inverse will vanish at the side faces according to the indicial roots.
The indicial roots for 
\begin{equation*}
	\Delta - s\lrpar{m-1-s}
\end{equation*}
are $s$ and $m-1-s$ and indeed we know from \cite{Mazzeo-Melrose} that its inverse
vanishes at the side faces to order $s$ (for this indicial root $x^s$ is locally in $L^2$ for $\Re(2s)>m-1$). 
The indicial operator for the Laplacian on forms is given in \cite[(3.2)]{Mazzeo:Hodge}. It preserves the splitting into tangential and normal parts of the form bundle and has indicial roots
\begin{equation*}
	r_{n} = \frac{m-1}{2} \pm \sqrt{\frac{(m-1)^2}{4} - k\lrpar{m-1-k} - \lambda},
	r_{t} = \frac{m-1}{2} \pm \sqrt{\frac{(m-1)^2}{4} - \lrpar{k-1}\lrpar{m-k} - \lambda},
\end{equation*}
on $k$ forms.
So if $\Re(s)$ is large enough, we will be able to apply Proposition \ref{Composition} and compose elements in the $even$ subcalculus. It turns out that this restricted composition is enough to show that the resolvent is in the even sub-calculus. We illustrate for the Laplacian on functions, but the same method extends to cover forms.

\begin{theorem}
If $g$ is an asymptotically hyperbolic metric which is $even$ mod $x^{2\ell}$ and $\Re(2s)> \max \lrpar{2\ell-m, m-1}$, then
\begin{equation*}
	\lrpar{\Delta - s\lrpar{m-1-s}}^{-1} \in \Psi^{2,\curly{E}}_{0,even}\lrpar{M;\Omega^{1/2}},
\end{equation*}
with $\Re(E_{11}) \geq 0$ and $\Re(E_{10}), \Re(E_{01}) \geq s$.
\end{theorem}

\begin{proof}
Let $A = \lrpar{\Delta - s\lrpar{m-1-s}}$ and assume that $\rho_{11}$ is a bdf for the front face in $C^{\infty}_{even}(M^2_0)$. 

We know that $N_{11}(A)^{-1}$ is an even operator on $\mathbb{H}^m$.
It extends to an $even$ operator
$C_{1} \in \Psi^{-2,\curly{H}^{(1)}}_{0,{even}}(M)$
with $\Re(H^{(1)}_{10}) \geq s$
satisfying
\begin{equation*}
	N^{0}_{11}(A)N^{0}_{11}(C_1) = Id.
\end{equation*}
Hence, $AC_1$ is in the $even$ subcalculus and
\begin{equation*}
	S_1 := AC_1 - Id 
	\in \rho_{11} \Psi^{0,\curly{G}^{(1)}}_{0,even}(M).
\end{equation*}
Iteratively, for $p<\ell$, assume we have found
\begin{equation*}
	\widetilde{C}_p = \sum_{j \leq p} \rho^{j-1} C_j, \phantom{x} 
	C_j \in \Psi^{-2,\curly{H}^{(j)} }_{0,{even}}(M)
\end{equation*}
 such that
\begin{equation*}
	S_p:= A\widetilde{C}_p - Id \in \rho_{11}^{p}
		\Psi_{0,even}^{0,\curly{G}^{(p)}}(M).
\end{equation*}
Then we find 
$C_{p+1} \in \rho_{11}^{2p}\phantom{x} \Psi^{-2,\curly{H}^{(j)} }_{0,{even}}(M)$
by extending $-N_{11}^0(A)^{-1}N^{2p}_{11}(S_p)$ off the front face in the $even$ subcalculus.
In which case
$\widetilde{C}_{p+1}= \widetilde{C}_p + C_{p+1}$ solves
\begin{equation*}
	S_{p+1}:= A\widetilde{C}_{p+1} - Id \in \rho_{11}^{p+1}
		\Psi_{0,{even}}^{0,\curly{G}^{(p+1)}}(M).
\end{equation*}

Thus we can find $C:=C_{2\ell}$ so that $S:=AC-Id$ vanishes to order $\rho_{11}^{2\ell}$ at the front face and
\begin{equation}\label{AInv}
	AC = Id + S \implies A^{-1} = C - A^{-1}S.
\end{equation}
Since $C$ is in the $even$ subcalculus, $A^{-1}$ will be $even$ if $A^{-1}S$ is.
The index set for $A^{-1}S$ at the front face is
\begin{equation*}
	(E_{11} + G^{(2\ell)}_{11}+2\ell) \bar{\cup} \lrpar{ E_{10} + G^{(2\ell)}_{01} + m},
\end{equation*}
and it follows from the construction that $\Re(G^{(p)}_{01}) \geq s$ for every $p$, so $A^{-1}S$ vanishes at the front face to order $2\ell$ and $A^-1$ is in the $even$ subcalculus.
\end{proof}

%%%%%%%%%%%%%%%%%%%%%%%%%%%%%%%%%%%%%%%%%%%%%%%%%%%%%%%%%%%%%%%%%%%%%%%%%%%%%%%%%%%%%%%%%%%%%%%%%%%%%%%%%%%%%%%%%%%%%%%%%%%%
\subsection{So is the Heat Kernel} $ $\newline
%%%%%%%%%%%%%%%%%%%%%%%%%%%%%%%%%%%%%%%%%%%%%%%%%%%%%%%%%%%%%%%%%%%%%%%%%%%%%%%%%%%%%%%%%%%%%%%%%%%%%%%%%%%%%%%%%%%%%%%%%%%%

Recall that the heat kernel of the Laplacian is, for any fixed $t>0$, an element of $\Psi_0^{-\infty}$.
Thus for fixed $t$, the heat kernel is given by a smooth function on the zero double space vanishing to infinite order on the side faces. We will use the result of the previous section to show that, for a metric $even$ mod $x^{2\ell}$,
the expansion of the heat kernel at the front face with respect to a special bdf has no odd terms below $x^{2\ell}$.

Choose $A, m \in \mathbb{R}^+$, and let $\gamma: \mathbb{R} \to \mathbb{C}$ be the path in the complex plane:
\begin{equation*}
      \gamma(t) = \begin{cases}
                             (-t, -A-mt) & \text{ if $t \leq 0$} \\
                             (t, -A-mt) & \text{if $t \geq 0$}
                            \end{cases}
\end{equation*}
Consider, for $t > 0$, the absolutely convergent integral:
\[ G(t) := \frac{1}{2\pi i} \int_{\gamma} e^{-t\lambda} \lrpar{\Delta - \lambda}^{-1} d\lambda. \]

As is well-known, this is the functional calculus expression for the heat kernel. It can be shown, as in \cite[(7.104)]{APS Book}, that $G(t)$ coincides with $e^{-t\Delta}$ as constructed in section $\mathcal{x}$\ref{sec:Construct}. Indeed, it is easy to see that, formally, for any smooth $f$,
\begin{equation*} \begin{split}
       \partial_t G(t)f &= %
              \frac{1}{2\pi i} \int_{\gamma} \partial_t e^{-t\lambda} \lrpar{\Delta - \lambda}^{-1}f d\lambda%
            =\frac{1}{2\pi i} \int_{\gamma} (-\lambda) e^{-t\lambda} \lrpar{\Delta - \lambda}^{-1}f d\lambda \\
         &=\frac{1}{2\pi i} \int_{\gamma} e^{-t\lambda} \lrpar{-\Delta+\lrpar{\Delta - \lambda}}%
              \lrpar{\Delta - \lambda}^{-1}f d\lambda %
           = -\Delta \lrpar{G(t)f} + \lrpar{ \frac{1}{2\pi i} \int_{\gamma} e^{-t\lambda} d\lambda } f \\
       & = -\Delta \lrpar{G(t)f}.
\end{split} \end{equation*}

This relates the integral kernels of the resolvent and heat kernel:
\begin{equation}\label{LapHeat}
       H_t(z,z') =
       \frac{1}{2\pi i} \int_{\gamma} e^{-t\lambda} R(\lambda, z, z') d\lambda .
 \end{equation}
As above, it is convenient to make the change of variables $\lambda = s(m-1-s)$, and to take $A, m$ large enough so that the resolvent is in the even calculus along the path $\tilde{\gamma}(t)$ (the image of $\gamma$ under $\lambda \mapsto s(m-1-s)$). Then the expansion of $H$ at the front face is obtained from the integral:
\begin{equation*}
       H_t(z,z') =
            \frac{1}{2\pi i} \int_{\tilde{\gamma}} e^{-ts(m-1-s)} R(s, z, z') (m-1-2s) ds. 
\end{equation*}
In particular, fix a bdf for the front face $\rho_{11} \in C^{\infty}_{odd}(M^2_0)$, then along this path $R(\tilde{\gamma}(t))$ will be $even$ and we conclude that the heat kernel itself will be even. We summarize with the following corollary.

\begin{corollary}\label{HeatKerFun}
For a conformally compact metric $even$ mod $x^{2\ell}$, the pointwise trace of the heat kernel is $even$ mod $x^{2\ell}$ for any special bdf, $x$.
In particular, if $2\ell \geq m$, the renormalized trace of the heat kernel on functions, obtained after choosing a special bdf, is independent of the choice of special bdf.
\end{corollary}

\begin{proof}
The results of the last section show that the heat kernel restricts to the diagonal to a density $even \mod x^{2\ell}$ (after identifying the diagonal of $M^2_0$ with $M$). As we discussed above, if $2\ell \geq m$ this implies that its renormalized integral is independent of the choice of special bdf (cf. \cite[Thm. 2.5]{My Preprint}).
\end{proof}

{\bf Remark.} One can play the same game \eqref{LapHeat} to define $f(\Delta)$ for other functions (with suitable growth restrictions, e.g., absolute integrability). These kernels inherit the even expansion of the resolvent at the front face. If their kernels are sufficiently regular to restrict to the diagonal, and $2\ell \geq m$ then their renormalized trace will be well-defined, independently of the choice of special bdf. A forthcoming article of Guillarmou explores this direction.

%%%%%%%%%%%%%%%%%%%%%%%%%%%%%%%%%%%%%%%%%%%%%%%%%%
%%%%%%%%%%%%%%%%%%%%%%%%%%%%%%%%%%%%%%%%%%%%%%%%%%
%%%%%%%%%%%%%%%%%%%%%%%%%%%%%%%%%%%%%%%%%%%%%%%%%%
%%%%%%%%%%%%%%%%%%%%%%%%%%%%%%%%%%%%%%%%%%%%%%%%%%
%%%%%%%%%%%%%%%%%%%%%%%%%%%%%%%%%%%%%%%%%%%%%%%%%%
%%%%%%%%%%%%%%%%%%%%%%%%%%%%%%%%%%%%%%%%%%%%%%%%%%
%%%%%%%%%%%%%%%%%%%%%%%%%%%%%%%%%%%%%%%%%%%%%%%%%%
%%%%%%%%%%%%%%%%%%%%%%%%%%%%%%%%%%%%%%%%%%%%%%%%%%
%%%%%%%%%%%%%%%%%%%%%%%%%%%%%%%%%%%%%%%%%%%%%%%%%%
%%%%%%%%%%%%%%%%%%%%%%%%%%%%%%%%%%%%%%%%%%%%%%%%%%
%%%%%%%%%%%%%%%%%%%%%%%%%%%%%%%%%%%%%%%%%%%%%%%%%%
%%%%%%%%%%%%%%%%%%%%%%%%%%%%%%%%%%%%%%%%%%%%%%%%%%

\end{document}